\documentclass[preprint, a4paper, 11pt]{elsarticle}
\usepackage[utf8]{inputenc}
\usepackage[official]{eurosym}

\usepackage{setspace}
\usepackage{amsmath,hhline,dcolumn,layout,booktabs, pdflscape}
\usepackage[margin=8pt,font=small,labelfont=bf]{caption}
\usepackage{algorithm}
\usepackage{algpseudocode}
\usepackage{floatflt}
\usepackage{graphicx}
\usepackage{subcaption}
\usepackage{longtable,verbatim}
\usepackage{enumitem,  rotfloat, supertabular}
\usepackage{lscape}
\usepackage{rotating}
\usepackage{booktabs}
\usepackage{amsfonts}
\usepackage{multirow}
\usepackage{chapterbib}
\usepackage{hyphenat}
\usepackage{multirow}
\usepackage{url}

\usepackage{geometry}
\newgeometry{vmargin={15mm}, hmargin={12mm,17mm}}

\usepackage[pdftex]{hyperref}
\hypersetup{colorlinks=true,allcolors=[rgb]{0,0,1}}

\usepackage{algorithm}
\usepackage{algpseudocode}

\usepackage{graphicx}

\usepackage{amsfonts,amscd}
\usepackage{amssymb,amsmath,amsthm}
\usepackage{mathtools}
\usepackage{bm, bbm}

\usepackage{natbib}
\bibpunct[, ]{(}{)}{;}{a}{,}{,}%
\def\bibsep{\smallskipamount}%
\title{An Optimization Framework for a Dynamic Multi-Skill Workforce Scheduling and Routing Problem with Time Windows and Synchronization Constraints}
\author[label1]{Onur Demiray}
\author[label2]{Doruk Tolga}
\author[label3]{Eda Y\"{u}cel\corref{cor1}}

\address[label1]{Department of Computing, Imperial College
London, Exhibition Rd, South Kensington, London SW7 2BX, United Kingdom. odemiray21@imperial.ac.uk}
\address[label2]{Department of Industrial Engineering, TOBB University of  Economics and Technology, S\"{o}ğ\"{u}t\"{o}z\"{u} Caddesi No:43, 06560, S\"{o}ğ\"{u}t\"{o}z\"{u}, Ankara, Turkey.  dtolga@etu.edu.tr}
\address[label3]{Department of Industrial Engineering, TOBB University of  Economics and Technology, S\"{o}ğ\"{u}t\"{o}z\"{u} Caddesi No:43, 06560, S\"{o}ğ\"{u}t\"{o}z\"{u}, Ankara, Turkey.  e.yucel@etu.edu.tr}
\cortext[cor1]{Corresponding Author: Eda Yücel,
e.yucel@etu.edu.tr}

\date{September, 2023}

\begin{document}

\begin{abstract}
This article addresses the dynamic multi-skill workforce scheduling and routing problem with time windows and synchronization constraints (DWSRP-TW-SC) inherent in the on-demand home services sector. In this problem, new service requests (tasks) emerge in real-time, necessitating a constant reevaluation of service team task plans. This reevaluation involves maintaining a portion of the plan unaltered, ensuring team-task compatibility, addressing task priorities, and managing synchronization when task demands exceed a team's capabilities.  To address the problem, we introduce a real-time optimization framework  triggered  upon the arrival of new tasks or the elapse of a set time. This framework redesigns the routes of teams with the goal of minimizing the cumulative weighted throughput time for all tasks. For the route redesign phase of this framework, we develop both a mathematical model and an Adaptive Large Neighborhood Search (ALNS) algorithm.  We conduct a comprehensive computational study to assess the performance of our proposed ALNS-based reoptimization framework and to examine the impact of reoptimization strategies,  frozen period lengths, and varying degrees of dynamism. Our contributions provide practical insights and solutions for effective dynamic workforce management in on-demand home services. 
\end{abstract}
\begin{keyword}
Dynamic vehicle routing, workforce scheduling and routing, reoptimization framework, adaptive large neighborhood search
\end{keyword}
\maketitle
\section{Introduction}
\label{section:introduction}
The increasing share of services is common for advanced economies as households tend to spend more on services as their income and wealth rise. Within the services sector there is a subsector called the on-demand home services sector that provides on-premise services primarily to households such as home health care or household professional services. The on-demand home services sector connects people who need assistance with day-to-day problems to those who are ready to provide instant solutions. According to Harvard Business Reports, these services attract 22.4 million consumers annually, who spend \$57.6 billion on them. The popular on-demand home services include health care, wellness, appliances, home cleaning, repairs and maintenance that require multi-skilled workforce. \citep{colby_bell_2017} 

The demand for home services increases with an increasing growth rate, and with the epidemic novel Covid-19, the adoption of the home services becomes extremely necessary not only to increase their benefits but also to provide relief to the people during the pandemic. Furthermore, home services assume a pivotal role in curbing and mitigating the transmission of this virulent virus. Their significance lies in the diminished need for extensive human interaction, which inherently aids in controlling the propagation of the virus.

On-demand home services can be either static or dynamic. In the static scenario, all tasks are known beforehand, while in the dynamic scenario, tasks show up dynamically as time progresses and the routes of the workforce should be continuously adapted to cater to evolving demand. The integration of mobile applications and diverse online platforms has led to the prevalence of  flexible, personalized, and responsive dynamic services. In practical terms, addressing dynamic challenges often involves treating them as sequences of static subproblems, as noted by \citet{CordeauLaporte2007}. The dynamism can be tackled through either offline (stochastic) or online (real-time) approaches in the solution methodology. In the offline scenario, robust solutions that consider dynamic demand data are produced before the operation. On the other hand, in the  online case, an initial solution is created based on the available data, which is then continually optimized as new service requests unfold during operations. This re-optimization from the current solution can be initiated whenever a new task arises. Alternatively, tasks can be buffered and periodically incorporated in the existing workforce routes in consolidated batches.

Effectively managing the operations of a dynamic on-demand on-site service system involves a multifaceted decision-making process encompassing task prioritization and identification of requisite skills, workforce dispatching, scheduling, and routing. Task prioritization is required due to the inherent variance in request importance and urgency. In the current landscape, service requests, or tasks, often demand distinct skill sets, rendering a service personnel incapable of addressing every task. This challenge is tackled through two different strategies: synchronization and team formation. Team formation finds applicability within static scenarios, if particularly when teams are established based on the demand data at the beginning of the planning phase. However, in the absence of synchronization, team formation alone falls short as a emergent request might necessitate skills beyond the capabilities of a single team (or worker) among the available team pool. Therefore, the subsequent problem of task assignment,  scheduling, and workforce routing   should account for synchronization to achieve a comprehensive solution.

The optimization problem addressed within this paper resides in the realm of dynamic workforce scheduling and routing problems with time windows and synchronization constraints (DWSRP-TW-SC), and its essence can be encapsulated as follows. In the static workforce scheduling and routing problem (WSRP), the aim is to allocate a collection of geographically dispersed service tasks with distinct skill prerequisites to a given ensemble of teams or workers, each possessing varying skill proficiencies. This entails identifying an optimal configuration of task-to-team assignments and devising a schedule for the task executions \citep{Cakirgil2020}. Transitioning to the dynamic rendition of this problem, the narrative takes on a more intricate dimension. New service requests arise concurrently with the execution of team task plans (TTPs). Considering the potential urgency of these new tasks, the TTP must be subject to continual re-evaluation and, if necessary, recalibration. In practice, when this recalibration takes place, a segment of the plan within a predefined timeframe - known as the frozen period - remains unaltered. In the course of re-optimizing the TTP, several factors merit attention. This includes the tasks situated within the frozen period (referred to as frozen tasks), the compatibility of team-task skills, task priorities, and the necessity for synchronization. All these elements converge to shape the trajectory of re-optimization.

To address the described DWSRP-TW-SC, we propose an optimization framework that is triggered whenever a predetermined number of new tasks arrive. The framework first identifies frozen tasks and determines the first available time and location of the team, and then re-optimizes the subsequent TTP with the objective of minimizing the total weighted throughput time of all tasks. For the route redesign phase of the framework, we develop both a mathematical model and a heuristic algorithm. 

The remainder of this paper is organized as follows. Section \ref{section:litreview} reviews the literature on related problems and approaches. The problem is formally described and mathematically modeled in Section \ref{section:problem}. Section \ref{section:solution} presents the proposed re-optimization framework that is based on an Adaptive Large Neighborhood Search (ALNS) heuristic. A comprehensive computational study is provided in Section \ref{section:compStudy}. Finally, Section \ref{section:conclusion} concludes the paper and draws directions for future research.

\section{Literature review}
\label{section:litreview}

In this section, we comprehensively examine the pertinent literature across three distinct categories: the studies focused on routing problems involving dynamic demand, the solution methodologies for dynamic routing problems, and the studies on measuring the level of dynamism.

\subsection{Routing problems with dynamic demand }\label{litreview_dynamicdemand}

The focus of this paper is the DWRSP-TW-SC problem, which is categorized as one of the dynamic demand vehicle routing problems (DVRP). As per \citet{Pillac2013}, routing problems with dynamic demand can be broadly classified into two primary categories: online and stochastic. In both cases, the demand for resources fluctuates over time. In the stochastic category, demand data is assumed to follow a probabilistic distribution, while in the online category, new demand occurrences are unpredictable. Our research falls into the dynamic and online category. Consequently, we will now provide an overview of the existing studies that belong to this specific classification.

\vspace{2mm}

The first study on DVRP that considered dynamic demand was \citet{WilsonColvin1977}. They developed a heuristic based on an insertion method for a single vehicle Dial-a-Ride problem, which concerns transportation services for disabled or elderly citizens. Since then, dynamic routing problems have gained traction in various fields, including food and human transportation. In recent years, there has been a notable surge in studies on dynamic routing problems, as noted by \citet{Pillac2013}. Notable reviews of DVRP and dynamic pickup and delivery studies have been provided by  \citet{Larsen2001}, \citet{Pillac2013},  and \citet{BERBEGLIA20108}. Many of these studies have focused on dynamic demand in Pickup and Delivery Problems (PDP) \citep{BERBEGLIA20108, Savelsbergh97, Gendreau2006}, city logistics \citep{Regnier-Coudert2016, Bieding2009, Campbell2005}, or human transportation issues like Dial-a-Ride Problems (DARP) \citep{CordeauLaporte2007_darp}. However, it's important to note that DARP studies that  consider dynamic demand \citep[e.g.,][]{Horn2002, Attanasio2004, COSLOVICH20061605, LOIS2017377} differ from our focus on DWSRP-TW-SC due to their homogeneous resources, such as vehicles or workforce. In some dynamic pickup and delivery studies \citep[e.g.,][]{Ferrucci2014}, heterogeneous vehicles are considered, but the demand structure, featuring origin and destination nodes, differs significantly from that of the DWSRP-TW-SC. 
\vspace{2mm}

Moving beyond these areas, as outlined in \citet{Pillac2013}, another application domain grappling with dynamic demand in vehicle routing problems is on-premise or home services. \citet{bostel2008multiperiod} explored dynamic technician routing and scheduling, considering a homogeneous workforce over a multi-period planning horizon. Their approach involved a memetic algorithm that schedules known tasks initially for each day and then seamlessly integrated new tasks into the existing plan. Additionally, they devised a column generation algorithm, but its applicability is limited to smaller-scale problems. Subsequently, \citet{Borenstein2008} addressed scheduling challenges with a multi-skilled workforce at British Telecom. They employed heuristic methods to assign technicians to task clusters and developed a rule-based system for task assignment when new tasks emerged. While their study tackled clustering and assignment aspects, it didn't delve into technician routing or task prioritization.

In recent times, \citet{Pillac2018} delved into the multi-skill technician scheduling and routing problem with dynamic demand, spare parts, and tool requirements. Their study treated tools as renewable resources and spare parts as non-renewable, consumed upon task completion. Technicians could replenish tools and spare parts at a central depot, offering them the flexibility to serve more requests. They introduced parallel processing, employing a parallel Adaptive Large Neighborhood Search (pALNS) algorithm to compute an initial solution and reoptimize it upon new request arrivals. What distinguishes \citet{Pillac2018} from our study is the consideration of spare part requirements and the capacity limitations of technicians regarding carrying spare parts. In our research, we assume sufficient resources within the vehicle and also allow synchronization, which involves tasks being serviced by multiple teams or technicians due to skill requirements.

\vspace{2mm}

\subsection{Solution approaches for dynamic routing problems}\label{litreview_solnapproaches}

Addressing routing problems in dynamic settings introduces an added layer of complexity, where response time emerges as a critical concern. Many of the solution methods typically employed for static problems prove to be excessively computationally intensive when applied to dynamic scenarios \citep{Gendreau2006}. 

\vspace{2mm}

Dynamic and deterministic problems are commonly tackled using two primary approaches: periodic reoptimization and continuous reoptimization \citep{Pillac2013}. Periodic reoptimization strategies involve generating vehicle routes at the outset of the planning horizon, often corresponding to a working day. As new data emerges throughout the day, or at predefined intervals, optimization techniques are deployed to update the routes based on this fresh information. For instance, \citet{Chang2003} explored real-time vehicle routing problems with time windows and simultaneous delivery/pickup demands, proposing a periodic reoptimization framework employing a Tabu Search (TS) algorithm. \citet{Montemanni2005} introduced an Ant Colony System algorithm for a real-world DVRP problem situated in the road network of Lugano city, based on data from a local fuel distribution company. \citet{Barcelo2007} developed a decision support system employing an iterative algorithm for dynamic routing of city logistic services. \citet{Attanasio2004} presented a parallel tabu-search heuristic within a periodic reoptimization framework to solve the dynamic multi-vehicle DARP. \citet{Pillac2018} also offered a periodic reoptimization framework, using a parallel ALNS, triggered each time a new request emerged in the context of dynamic technician routing and scheduling.

\vspace{2mm}
In contrast, continuous reoptimization approaches run continuously throughout the day, relying on adaptive memory \citep{Taillard2001} to store alternative solutions for adapting to new data without complete reoptimization. \citet{Gendreau1999} developed a parallel TS with adaptive memory for a dynamic VRP with time windows. This approach adapted the static TS to dynamic scenarios through parallel execution, where threads running in the background are interrupted upon data updates (e.g., new requests or completed customer services). Subsequently, the solutions in the adaptive memory are updated, and the parallel search process restarts with this updated information. This framework was extended to dynamic VRP by \citet{Ichoua2000} and \citet{Ichoua2003}. \citet{Bent2004} explored a dynamic VRP with time windows and stochastic customer information, introducing a multiple plan approach (MPA) inspired by \citet{Gendreau1999}'s framework. The MPA maintains a pool of routing plans consistent with current decisions and removes incompatible ones from the pool. \citet{Benyahia1998} proposed a Genetic Algorithm (GA) modeling the vehicle dispatcher's decision-making process for dynamic PDP. GA techniques are employed to find a utility function that approximates the dispatcher's choices. \citet{Hemert2004}, \citet{HAGHANI2005}, and \citet{Cheung2008} extended the continuous reoptimization framework using GAs for dynamic PDPs. These GAs, designed for dynamic contexts, closely resemble those used for static problems but operate continuously throughout the planning horizon, adapting solutions as input data changes \citep{Pillac2013}.

\vspace{2mm}
In this study, we propose a periodic reoptimization approach based on an ALNS. Our optimization framework is triggered when a predetermined number of new tasks, denoted as $\beta^{\text{task}}$, arrive, or when a set amount of time, referred to as $\beta^{\text{time}}$, elapses. We conduct a computational analysis to assess various values of $\beta^{\text{task}}$ and $\beta^{\text{time}}$ (see Section~\ref{section:compStudy_reoptimizationPeriodComparison}). To account for teams engaged in task execution or en route to task locations at the time of reoptimization, our proposed framework freezes team-task assignments within a designated frozen period, denoted as $f$, starting from the reoptimization moment. We also conduct experiments to evaluate the optimal length of the frozen period (see Section~\ref{section:compStudy_frozenPeriodComparison}).

\subsection{Measuring dynamism}\label{litreview_measuringDynamism}

Problems with dynamic data, even within the same problem class, can exhibit varying degrees of dynamism. According to \citet{Ichoua2007}, the dynamism of a problem instance can be characterized using two key dimensions: the frequency of changes and the urgency of customer requests. The first dimension relates to how quickly new information emerges, while the second pertains to the time interval between the arrival of a new request and its expected service time or deadline. Based on these considerations, three distinct metrics have been proposed to measure the dynamism of a problem or instance: $(i).$ Degree of Dynamism ($\delta$), as introduced by \citet{Lund1996}, is defined as the ratio of new requests ($n_{\text{d}}$) to the total number of requests ($n_{\text{tot}}$), i.e., $\delta = n_{\text{d}}/n_{\text{tot}}$.
$(ii).$ Effective Degree of Dynamism ($\delta^e$), as proposed by \citet{Larsen2001}, takes into account the arrival times of requests ($a_i$) and the length of the planning horizon ($\tau^{\max}$). For requests known in advance, the arrival time is considered as 0. $\delta^e$ is calculated as follows: $\delta^e = \frac{1}{n_{\text{tot}}}\sum_{i\in N} \frac{a_i}{\tau^{\max}}$, where $N$ represents the set of requests.
$(iii).$ Effective Degree of Dynamism with Time Windows ($\delta^e_{\text{tw}}$), also proposed by \citet{Larsen2001}, takes into consideration both the arrival time and the latest start time of requests ($a_i$ and $l_i$, respectively), along with the length of the planning horizon ($\tau^{\max}$). It is calculated as: $\delta^e_{\text{tw}} = \frac{1}{n_{\text{tot}}}\sum_{i\in N} \left(1 - \frac{(l_i - a_i)}{\tau^{\max}}\right)$.

\citet{Larsen2001} further classified dynamic problems or instances into three categories based on the effective degree of dynamism ($\delta^e$): Weekly Dynamic Systems having $\delta^e \leq 0.3$,
Moderately Dynamic Systems having $0.3 < \delta^e \leq 0.8$,
Strongly Dynamic Systems where $\delta^e > 0.8$.
In our study, we investigate the impact of the degree of dynamism on the quality of solutions provided by our proposed reoptimization framework through computational analysis (see Section~\ref{section:compStudy_largeSizeDWSRPinstances}).

\vspace{2mm}

\section{Problem description and mathematical model}
\label{section:problem}

DWRSP-TW-SC entails the allocation of crews equipped with certain skills to a geographically-dispersed set of tasks that are randomly revealed during the day. The problem involves a predetermined and unalterable set of crews denoted by $\mathcal{K}$, corresponding to the technician compositions, therefore the skill qualifications of crews are given at the beginning of the day and remain constant throughout the course of the planning horizon. On the other hand, the tasks are dynamic, and new tasks with varying skill requirements and priorities emerge as the operation of the crews progresses in the field. The sequence of tasks to be completed by each crew is defined by the so-called Team Task Plan (TTP), which must be updated whenever a predetermined number of new tasks appear. This process of updating the TTP is referred to as a \textit{re-optimization  problem (period/phase)}, which is regarded as a separate optimization problem.

\vspace{2mm}

The $t^{\textit{th}}$ re-optimization period begins at time $\tau^{(t)}$. At this time point, the system's tasks can be categorized into two sets: (i) $\mathcal{N}_{prev}^{(t)}$, corresponding to the tasks that have not been addressed yet and were available at the previous re-optimization period, and (ii) $\mathcal{N}_{new}^{(t)}$, corresponding to the tasks that have accumulated in the system between the previous and current re-optimization periods. In other words, $\mathcal{N}_{new}^{(t)}$ is made up of tasks whose arrival times fall within the time interval $(\tau^{(t-1)}, \tau^{(t)}]$. At this point, it is worth noting that the initial static optimization problem at the start of the day can be considered as the previous re-optimization period for the first re-optimization period. Hence, without loss of generality, we refer to this initial static problem as the $0^{\textit{th}}$ re-optimization problem that occurs at the time $\tau^{(0)}=0$.

\vspace{2mm}

The team assignment of the tasks in $\mathcal{N}_{prev}^{(t)}$ can be changed except a subset of tasks in $\mathcal{N}_{prev}^{(t)}$, which are planned to be started by a crew within a specified time period called \textit{frozen period}. The frozen period of an appropriate duration aims to prevent confusions in switching from one TTP to another TTP. Let $f^{(t)}$ denote the duration of the frozen period for the $t^{\textit{th}}$ re-optimization period. Then, $[\tau^{(t)}, \tau^{(t)} + f^{(t)}]$ specifies the frozen period of the $t^{\textit{th}}$ re-optimization period. Tasks whose start times fall into this period are called \textit{frozen tasks} and are referred to as $\mathcal{N}_{frozen}^{(t)}$.  During the $t^{\textit{th}}$ re-optimization period, it is possible for the completion time of the final frozen task to differ among crews. As a result, the starting time ($\eta_k^{(t)}$ for crew $k \in \mathcal{K}$) of each crew to the same period may vary. The starting location ($\phi_k^{(t)}$) of crew $k \in \mathcal{K}$ is determined by the location of its final frozen task. In other words, the last task completed by a crew during the frozen period serves as the basis for determining its starting location to the corresponding re-optimization period.

\vspace{2mm}

The problem to be solved for the $t^{\textit{th}}$ re-optimization period is defined on a directed graph $G^{(t)}= (\mathcal{V}^{(t)}, \mathcal{A}^{(t)})$, where $\mathcal{V}^{(t)} = \mathcal{N}^{(t)}  \cup  \Phi \cup \{0\} $ denotes  the set of vertices (nodes), $\mathcal{N}^{(t)}$ 
corresponds to the whole set of tasks that should be considered at the time of re-optimization, i.e., $\mathcal{N}^{(t)} = \mathcal{N}_{prev}^{(t)}\setminus \mathcal{N}^{(t)}_{frozen} \cup \mathcal{N}_{new}^{(t)}$, $\Phi$ corresponds to the starting locations of crews, i.e., $\Phi = \bigcup_{k\in \mathcal{K}} \phi_k$, node $0$ represents the central depot that each team  should go at the end of the work day, and $\mathcal{A}^{(t)}=\{(i,j):i \in \mathcal{V}^{(t)}, j \in \mathcal{V}^{(t)} \setminus \Phi, i \neq j\}$ denotes the set of arcs between the nodes.

\vspace{2mm}

For each arc $(i,j) \in A^{(t)}$, $c_{ij}$ denotes the traveling time from node $i\in \mathcal{V}^{(t)}$ to $j\in \mathcal{V}^{(t)}\setminus\Phi$. Each task $i \in \mathcal{N}^{(t)}$ is associated with an arrival time $a_i$, processing time $p_i$,  priority $w_{i}$, and time window $[e_i, l_i]$. The set of skills is denoted by $\mathcal{Q}$. The $0-1$ parameter $u_{iq}$ indicates whether task $i \in \mathcal{N}^{(t)}$ requires skill $q \in \mathcal{Q}$ and $v_{kq}$ indicates whether crew $k \in \mathcal{K}$ has skill $q \in \mathcal{Q}$. There may be some tasks that require a skill set that cannot be satisfied by a single team. Depending on their skill requirements, two or more crews are assigned to these task and these crews have to be synchronized, meaning that they need to start at the same time. 

\vspace{2mm}

The planning horizon, denoted by $[\tau^{(0)}=0, \tau^{\max}]$, spans from the beginning of the working day until its conclusion at $\tau^{\max}$. If a task cannot be completed within this timeframe, it is outsourced and assumed to be accomplished by the end of the day, i.e., at $\tau^{\max}$. Our objective is to minimize the total weighted throughput time of the tasks, where the throughput time of a task is defined as the duration between its arrival and completion. In the remainder of this section, we will initially present an illustrative example to demonstrate how dynamic tasks are executed by re-optimization problems during the day, followed by the introduction of a mixed-integer programming formulation designed to solve any given re-optimization problem.

\vspace{2mm}

\subsection{An illustrative example}\label{illustrativeExample}

\vspace{2mm}

Within this section, we present an illustrative example involving two crews, eight tasks, and three skills, with the primary objective of showcasing the execution of dynamically revealed tasks through different re-optimization phases. By examining this scenario, we aim to provide a comprehensive understanding of how our dynamic framework effectively handles the allocation and sequencing of tasks that unfold over time. Additionally, we aim to demonstrate how synchronization due to the assignment of multiple teams to a task might cause idle times for crews.

\vspace{2mm}

\begin{figure}[htbp]
  \centering
  \includegraphics[width=0.7\linewidth]{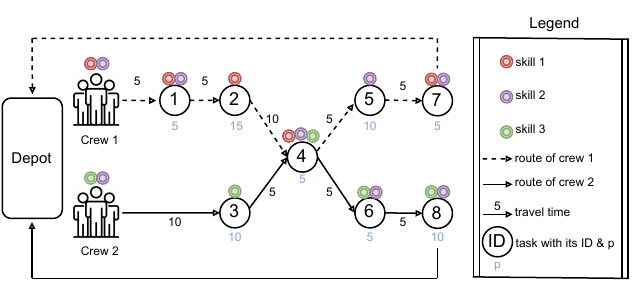}
  \caption{An illustrative solution for an instance having two crews, 8 tasks, and 3 skills}
  \label{fig:illustrative_example}
\end{figure}

The hypothetical example is visually represented in Figure \ref{fig:illustrative_example}, showing the routes taken by the crews. In this example, crew 1 executes tasks in the sequence of $ <1,2,4,5,7>$, while crew 2 follows the route of $<3,4,6,8>$. The skills are represented by colored double line circles positioned above the crew and task symbols. Specifically, there are three distinct skills denoted by different colors: red, purple, and green. Figure \ref{fig:illustrative_example} provides insight into the skill composition of the crews. Crew 1 possesses skills associated with the red and purple circles, enabling them to complete tasks 1, 2, 5, and 7 autonomously. Conversely, crew 2 possesses skills indicated by the purple and green circles, equipping them to undertake tasks 3, 6, and 8 independently. However, task 4 necessitates collaboration between both crews since it requires the combined skill set. Consequently, the crew that reaches the location of task 4 earlier must await the arrival of the other crew, highlighting the concept of synchronization.

\vspace{2mm}

\begin{figure}[htbp]
  \centering
  \includegraphics[width=0.6\linewidth]{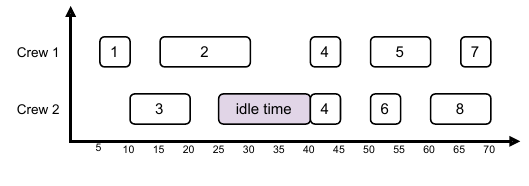}
  \caption{Gantt chart of the crews for the solution of the illustrative example represented in Figure \ref{fig:illustrative_example}}
  \label{fig:illustrative_example_gantt}
\end{figure}

In Figure \ref{fig:illustrative_example}, the numbers displayed over the arcs indicate the corresponding travel times required to traverse those arcs. Additionally, the blue-colored numbers beneath the nodes represent the process time of each respective task. Consequently, the Gantt chart in Figure \ref{fig:illustrative_example_gantt} visualizes the start and finish times of each task. For tasks completed by a single team, the start time can be easily determined by adding the completion time of the previous task to the travel time required to reach the current task when all tasks are already revealed. However, in cases where synchronization is required, such as in Task 4, the team that arrives earlier experiences idle time. For instance, in order to complete Task 4, crew 2 arrives at time 25 while crew 1 arrives at 40. Consequently, crew 1 has an idle time of 15 time units while waiting for crew 2 to arrive.

For illustrative purposes, let us assume that the given routes and schedules shown in Figure \ref{fig:illustrative_example} pertain to the static problem. While we here demonstrate the transition to the first re-optimization, the same approach can be applied to subsequent transitions between consecutive re-optimization periods. At the initial time $0$, the set $\mathcal{N}^{(0)}$ comprises tasks $\{1,2,3,4,5,6,7,8\}$. Assuming $\tau^{(1)}=45$ and $f^{(1)}=10$, the frozen period spans the interval $[45,55]$. Prior to reaching $\tau^{(1)}=45$, tasks 1, 2, 3, and 4 have already been completed, resulting in $\mathcal{N}_{prev}^{(1)}$ being defined as $\{5,6,7,8\}$. During the frozen period, crew 1 initiates task 5 and completes it at time 60. Consequently, crew 1 commences the first re-optimization period from the location of task 5 at time 60, i.e., $\phi_1^{(1)}=5$ and $\eta_1^{(1)}=60$. Employing the same logic, we have $\phi_2^{(1)}=6$ and $\eta_2^{(1)}=55$. Therefore, $\mathcal{N}_{frozen}^{(1)}$ is defined as $\{5,6\}$. Suppose that only two tasks, task 9 and task 10, appear within the time interval $(0,45]$. Consequently, $\mathcal{N}^{(1)}$ is defined as $\{7,8,9,10\}$. These four tasks must be addressed during the first re-optimization phase.

\vspace{2mm}
\subsection{Mixed integer programming formulation}\label{mip}
\vspace{2mm}

In this section, we present the mixed-integer programming formulation for the $t^{\textit{th}}$ re-optimization problem of DWRSP-TW-SC for any $t$. We start with introducing the  decision variables:

\vspace{-4mm}

\begin{align*}
    X_{ij}^k &:
    \begin{cases}
    1, & \text{if arc }(i,j)\in\mathcal{A}^{(t)} \text{ is traversed by crew }k\in\mathcal{K}\\
    0, & \text{otherwise}
    \end{cases} \\
    Y_i^k &:
    \begin{cases}
    1, & \text{if node }i \in \mathcal{V}^{(t)}\setminus\{0\} \text{ is assigned to crew }k\in\mathcal{K}\\
    0, & \text{otherwise}
    \end{cases} \\
    O_i &: 
    \begin{cases}
    1, & \text{if task }i \in \mathcal{N}^{(t)} \text{ is accomplished via outsourcing}\\
    0, & \text{otherwise}  
    \end{cases}\\
    D_{i}^k &: \text{departure time of crew }k \in \mathcal{K} \text{ from node }i\in\mathcal{V}^{(t)}\\
    C_i &: \text{completion time of task } i \in \mathcal{N}^{(t)} 
\end{align*}

 Then, the $t^{\textit{th}}$ re-optimization problem at any time in between $\tau^{(0)}$ and $\tau^{\max}$ has a planning horizon of $[\tau^{(t)}, \tau^{\max}]$ and  can be formulated as follows:

\begin{align}
        \text{minimize} & \quad \displaystyle \text{TWTT} =  \sum_{i \in \mathcal{N}^{(t)}} w_i (C_i - a_i) \label{mip_objective}\\
        \mathrm{s.t.} & \displaystyle\quad O_i + \sum_{k \in \mathcal{K}}Y_{i}^k \geq 1 & \forall i \in \mathcal{N}^{(t)} \label{c1:job_must_be_completed}\\
        & \displaystyle \quad  \sum_{k \in \mathcal{K}}Y_{i}^k \leq |\mathcal{K}| (1-O_i) & \forall i \in \mathcal{N}^{(t)} \label{c2:if_outsourced_no_team} \\
        & \displaystyle \quad Y_{\phi_k^{(t)}}^k = 1 & \forall k \in \mathcal{K} 
        \label{c3:start_location} \\
        & \displaystyle \quad D_{\phi_k^{(t)}}^k = \eta_k^{(t)} & \forall k \in \mathcal{K} \label{c4:start_time}\\
        & \displaystyle \quad \sum_{i \in \mathcal{N}^{(t)} \cup \{\phi_k\}}X_{i0}^k = 1 & \forall k \in \mathcal{K} 
        \label{c5:return_depot} \\
        & \displaystyle \quad \sum_{\substack{j \in \mathcal{V}^{(t)}:\\(i,j) \in \mathcal{A}^{(t)}}} X_{ij}^k = Y_i^k & \forall i \in \mathcal{N}^{(t)} \cup \{\phi_k\}, k \in \mathcal{K} \
        \label{c6:X_Y_match} \\
        & \displaystyle \quad \sum_{\substack{j \in \mathcal{N}^{(t)} \cup \{\phi_k\}:\\(j,i) \in \mathcal{A}^{(t)}}}X_{ji}^k - \sum_{\substack{j \in \mathcal{N}^{(t)} \cup \{0\}:\\(i,j) \in \mathcal{A}^{(t)}}}X_{ij}^k = 0 & \forall i \in \mathcal{N}^{(t)}, k \in \mathcal{K} \label{c7:flow_balance} \\
        & \displaystyle \quad \sum_{k \in \mathcal{K}} v_{kq}Y_{i}^k \geq u_{iq}- O_i & \forall i \in \mathcal{N}^{(t)}, q \in \mathcal{Q} \label{c8:skill_competence}\\
        & \displaystyle \quad D_i^k + c_{ij} + p_j \leq D_j^k + \tau^{\max}\color{black}(1-X_{ij}^k) &\forall k \in \mathcal{K}, i \in \mathcal{V}^{(t)}\setminus\{0\}, j \in \mathcal{V}^{(t)}\setminus\{\phi_k\}: i \neq j \label{c9:departure_time_of_crews} \\
        & \displaystyle \quad D_i^k + \tau^{\max}(1-Y_i^k) \geq C_i - \tau^{\max}O_i & \forall i \in \mathcal{N}^{(t)}, k \in \mathcal{K} \label{c10:task_completion_time_1}\\
        & \displaystyle \quad D_i^k \leq C_i & \forall i \in \mathcal{N}^{(t)}, k \in \mathcal{K} \label{c11:task_completion_time_2}\\
        & \displaystyle \quad C_i \leq \tau^{\max}O_i + (l_i+p_i)(1-O_i) \quad & \forall i \in \mathcal{N}^{(t)} \label{c12:time_window_1}\\
        & \displaystyle \quad C_i \geq \tau^{\max}O_i + (e_i+p_i)(1-O_i) \quad & \forall i \in \mathcal{N}^{(t)} \label{c13:time_window_2}
            \end{align}
        \begin{align}
        & \displaystyle \quad X_{ij}^k \in \{0,1\} & \forall (i,j) \in \mathcal{A}^{t}, k \in \mathcal{K} \label{c14:domain_X} \\
        & \displaystyle \quad Y_{i}^k \in \{0,1\} & \forall i \in \mathcal{V}^{t} \setminus\{0\}, k \in \mathcal{K} \label{c15:domain_Y} \\
        & \displaystyle \quad O_i \in \{0,1\} & \forall i \in \mathcal{N}^{t}\label{c16:domain_O} \\
        & \displaystyle \quad D_{i}^k \geq 0 & \forall i \in \mathcal{V}^{t} \setminus\{0\}, k \in \mathcal{K} \label{c17:domain_D} \\
        & \displaystyle \quad C_i \geq 0 & \forall i \in \mathcal{N}^{t}\label{c18:domain_C} 
    \end{align}
 
The objective function (\ref{mip_objective}) minimizes the weighted sum of task throughput times for all tasks. It's important to note that when we remove the constant term from the objective, it essentially transforms into the minimization of the weighted sum of completion times for all tasks. As such, these two descriptions of the objective are used interchangeably throughout this article.  Constraints (\ref{c1:job_must_be_completed}) ensure that all tasks in the system are either outsourced or assigned to at least one crew. Depending on the capabilities of the crews, multiple crews may be assigned to a task. However, if a task is outsourced, no crew can be assigned to it, as stated in constraints (\ref{c2:if_outsourced_no_team}). Each crew begins their tour from their last completed task, as determined by constraints (\ref{c3:start_location}) and (\ref{c4:start_time}), and must complete their tour at the central depot, as indicated by constraint (\ref{c5:return_depot}). Constraints (\ref{c6:X_Y_match}) imply that if a node is assigned to a crew, the task must be part of that crew's route. Constraints (\ref{c7:flow_balance}) ensure flow conservation, ensuring that a crew arrives at and departs from the node of an assigned task. Constraints (\ref{c8:skill_competence}) guarantee that the skill requirements of in-house tasks are satisfied by the assigned team(s). Constraints (\ref{c9:departure_time_of_crews}) calculate the departure time of the crews from the nodes. Constraints (\ref{c10:task_completion_time_1}) and (\ref{c11:task_completion_time_2}) determine the completion time of in-house tasks. Constraints (\ref{c12:time_window_1}) and (\ref{c13:time_window_2}) set the completion time of outsourced tasks to the end of the day and enforce time window constraints for non-outsourced tasks. Finally, Constraints (\ref{c14:domain_X}) - (\ref{c18:domain_C}) specify the domains of the decision variables.

\section{Solution approach}
\label{section:solution}

Considering the inherent nature of the problem, it is imperative that the $t^{th}$ re-optimization phase, which commences at time $\tau^{(t)}$, is concluded by time $\tau^{(t)} + f^{(t)}$. While the mathematical model can produce optimal solutions for small- to medium-sized problem instances within this time window, obtaining the optimal solution through the model becomes increasingly time-intensive for larger instances. As a result, we propose an ALNS heuristic to generate high-quality solutions within the designated time interval. The proposed ALNS algorithm initiates with an initial solution constructed by a constructive heuristic, the details of which are presented in Section \ref{subsec:constructive}. The subsequent Section \ref{subsec:alns} provides a detailed description of the ALNS algorithm.

\subsection{Constructive heuristic}\label{subsec:constructive}

The algorithm starts with sorting the tasks in $\mathcal{N}^{(t)}$ in a non-increasing order based on their priorities. Considering that a single team may not be able to fulfill a task's requirements, there can exist multiple team combinations that are capable of satisfying the skill requirements of the task. Subsequently, starting from the first task in the ordered list, the algorithm identifies the \textit{irreducible competible team combinations} for each task, taking into account its skill requirements. These combinations represent the alternative minimal sets of teams capable of fulfilling the skill requirements of the task. The collection of all those combinations for task $i$ forms the set $\mathcal{A}_i$.

\vspace{2mm}

For instance, consider a scenario where there are five distinct skills, and the task being considered, task $i$, necessitates the first, second, and fifth skills, represented as $u_i = [1, 1, 0, 0, 1]$. Assuming the existence of three teams with corresponding skill qualifications, denoted as $v_1 = [1, 1, 1, 0, 1]$, $v_2 = [1, 0, 0, 0, 0]$, and $v_3 = [0, 1, 0, 0, 1]$, the set of irreducible competible team combinations $\mathcal{A}_i$ is derived as $\{ \{1\}, \{2, 3\}\}$.

\vspace{2mm}

If no feasible team combination is identified for a given task $i$, denoted as $\mathcal{A}_i = \emptyset$, the task is subsequently outsourced. Next, for each available alternative team combination, the optimal insertion point is determined for each team within the combination, taking into consideration their respective time windows and the resulting total insertion cost. Following this, the task is inserted at the best insertion point within each team of the selected team combination, which exhibits the lowest total insertion cost. This iterative process continues until all tasks have been assigned to teams or outsourced. For a comprehensive understanding of the constructive heuristic's implementation, refer to Algorithm \ref{Algorithm:CH}, which provides the pseudocode outlining the step-by-step procedure.

\begin{algorithm}[ht!]
\caption{Constructive Heuristic}
\begin{algorithmic}[1]
\Require  {$\mathcal{N}^{(t)}$: set of tasks to be considered at the time of re-optimization,  $\mathcal{K}$: set of crews}
\State $\mathcal{N}^{outsourced} \leftarrow \emptyset$
\State $\mathcal{N}^{sorted} \leftarrow SortTasks(\mathcal{N}^{(t)})$
\For {each $i \in \mathcal{N}^{sorted}$}
   \State $\mathcal{A}_i \leftarrow FindIrreducibleCompetibleCombinations(i,\mathcal{K})$
   \If {$\mathcal{A}_i = \emptyset$}
   \State $\mathcal{N}^{outsourced} \leftarrow \mathcal{N}^{outsourced}  \cup$ $\{i\}$
   \Else
   \State ${\mathcal{A}_i}^* \leftarrow FindBestCombination(\mathcal{A}_i)$
\EndIf
\For {each $k \in \mathcal{A}_i^*$}
     \State Insert task $i$ to the best insertion point of the route of crew $k$
\EndFor
\EndFor
\Ensure $x^0$: initial solution\end{algorithmic}\label{Algorithm:CH}
\end{algorithm}

\subsection{Adaptive large neighborhood search algorithm}\label{subsec:alns}

ALNS was originally introduced by Ropke and Pisinger in their works \cite{RopkePisinger2006a, RopkePisinger2006b}. This metaheuristic approach operates through iterative selection of a destroy and repair heuristic from a designated set, based on their past performance. By employing this strategy, the ALNS generates candidate solutions and determines their acceptance or rejection through a probabilistic criterion. Renowned for its ability to balance exploration and exploitation, the ALNS has demonstrated its effectiveness in optimizing various problem domains, particularly in the realm of routing problems. Thus, we propose the utilization of an ALNS algorithm to enhance the initial solution derived from the constructive heuristic, leveraging its strengths to further refine the solution quality.

\vspace{2mm}

The pseudocode for our Adaptive Large Neighborhood Search algorithm is presented in Algorithm \ref{Algorithm:ALNS}. The algorithm begins by constructing an initial solution $x^0$ using the constructive heuristic. It then initializes the current and best-found solutions, while assigning equal weights to the destroy and repair heuristics. In each iteration, the algorithm selects a destroy heuristic $h^{\text{destroy}}$ and a repair heuristic $h^{\text{repair}}$ based on their performance in previous iterations, taking their weights into account. Given the current solution $x^{\text{current}}$, the selected destroy heuristic removes a certain number of tasks from the assigned teams, and the repair heuristic reinserts the removed tasks into the routes, if possible. If the resulting solution $x^{\text{candidate}}$ satisfies the stochastic acceptance criteria, it replaces the current solution $x^{\text{current}}$. Furthermore, if $x^{\text{candidate}}$ is superior to the best solution found so far, denoted as $x^{\text{best}}$, $x^{\text{candidate}}$ replaces $x^{\text{best}}$ as the new best solution. The algorithm then updates the weights of the destroy and repair heuristics based on the performance of $h^{\text{destroy}}$ and $h^{\text{repair}}$, and proceeds to the next iteration. The weights of the destroy and repair heuristics influence their selection probability in subsequent iterations. The algorithm continues until the stopping criterion is met.

\begin{algorithm}[ht!]
\caption{ALNS}
\begin{algorithmic}[1]
\Require $x^0$: initial solution 
\State $x^{\text{best}}\leftarrow x^0$
\State $x^{\text{current}}\leftarrow x^0$
\State $\mathbf{w^{\text{destroy}}} \leftarrow \text{InitializeDestroyWeights()}$
\State $\mathbf{w^{\text{repair}}} \leftarrow \text{InitializeRepairWeights()}$
\State \textit{iter} $\leftarrow 1$
\While {\text{stopping criteria is not met}}
\State $h^{\text{destroy}} \leftarrow \text{ChooseDestroyHeuristic}(\mathbf{w^{\text{destroy}}})$ 
\State $h^{\text{repair}} \leftarrow \text{ChooseRepairHeuristic}(\mathbf{w^{\text{repair}}})$ 
\State $x^{\text{candidate}} \leftarrow h^{\text{repair}}(h^{\text{destroy}}(x^{\text{current}}))$
\If{Accept$(x^{\text{candidate}},x^{\text{current}} )$}
    \State $x^{\text{current}} \leftarrow x^{\text{candidate}}$
    \If{$\text{TWTT}(x^{\text{candidate}})<\text{TWTT}(x^{\text{best}})$}
    \State $x^{\text{best}} \leftarrow x^{\text{candidate}}$ 
    \EndIf
\EndIf
\State $\mathbf{w^{\text{destroy}}} \leftarrow \text{UpdateDestroyWeights()}$
\State $\mathbf{w^{\text{repair}}} \leftarrow \text{UpdateRepairWeights()}$
\State $iter \leftarrow iter+1$
\EndWhile
\Ensure $x^{\text{best}}$: best found solution
\end{algorithmic}\label{Algorithm:ALNS}
\end{algorithm}

The remainder of this section is dedicated to describing the components of the ALNS algorithm, including the destroy and repair heuristics, the acceptance criteria, the heuristic weight adjustment and selection processes.

\subsection{Destroy heuristics}

Our ALNS incorporates five destroy heuristics, denoted as the set $\mathcal{D}$. These heuristics include problem-specific variations of the random, worst, and related removal heuristics proposed by \citet{RopkePisinger2006b},  \citet{PisingerRopke2007}, \citet{shaw1997}, and \citet{shaw1998}. The destroy heuristics employed in our algorithm are random task removal, worst task removal, random team removal, worst team removal, and Shaw removal. Prior to the application of a destroy heuristic, we determine the \textit{degree of destruction}, denoted as $\gamma$. This parameter specifies the percentage of tasks to be removed from the current solution $x^{\text{current}}$. At each iteration, $\gamma$ is assigned a random value between $\gamma^{min}$ and $\gamma^{max}$.

\subsubsection{Random task removal (RTaR):} \label{subsubsec:rtr} The algorithm employs a randomized selection mechanism to choose a task, which is subsequently removed from the routes of all teams to which it is assigned. This iterative process continues until $\gamma$ percentage of the total tasks has been successfully eliminated from the current solution.

\subsubsection{Worst task removal (WTaR):} The algorithm computes a score for each task based on its impact on the objective function. Through this scoring process, the task with the highest score is identified as the worst task, as we have a minimization problem. Subsequently, the algorithm removes this worst task from the routes of all teams to which it is assigned. The procedure continues in an iterative manner until $\gamma$ percentage of the total tasks has been successfully eliminated from the current solution.

\subsubsection{Random team removal (RTeR):} \label{subsubsec:rter} The algorithm randomly chooses a team and proceeds to remove all tasks from its route until $\gamma$ percentage of the total tasks has been eliminated. It is important to note that if a task that needs to be removed is assigned to multiple teams, it will be removed from all teams to which it is assigned.

\subsubsection{Worst team removal (WTeR):} 
The algorithm evaluates the impact of each team on the objective function and assigns a score accordingly. Through this scoring process, the team with the highest score is identified as the worst team. The algorithm then proceeds to remove the tasks associated with the identified worst team from the routes of all teams to which they are assigned. This iterative process continues until $\gamma$ percentage of the total tasks has been successfully eliminated. 

\subsubsection{Shaw removal (SR): } 

According to the Shaw Removal heuristic, it is more promising to remove tasks that are related rather than unrelated, as it is more likely to lead to an improved solution. The relatedness between two tasks, $R(i, j)$, is calculated using the $L_2$ norm distance between the skill vectors of the tasks, as described by Equation \ref{Rij}.

\begin{align}
R(i, j) = \sqrt{\sum_{q\in \mathcal{Q}} ((u_{iq})^2 -(u_{jq})^2)} \label{Rij}
\end{align}

As the value of $R(i, j)$ approaches 0, the tasks $i \in \mathcal{N}^{(t)}$ and $j \in \mathcal{N}^{(t)} \setminus \{i\}$ are considered more closely related. The Shaw Removal (SR) heuristic initiates by randomly selecting a task $i$ from the set of tasks $\mathcal{N}^{(t)}$. It then removes task $i$ and identifies the most related task $j^*$, where $j^* = \arg\min_{j \in \mathcal{N}^{(t)} \setminus \{i\}} R(i,j)$.  This iterative process continues until $\gamma$ percentage of the total tasks has been successfully eliminated. 

\subsection{Repair heuristics}
\label{sec:insertion}

The proposed ALNS algorithm incorporates three distinct repair heuristics, represented by the set $\mathcal{R}$: random insertion, greedy insertion, and regret insertion.

\vspace{2mm}

Let $\mathcal{N}^{(t)}_{or}$ denote the amalgamation of outsourced tasks and the set of removed tasks resulting from the application of the selected removal heuristic. Prior to executing a repair heuristic, the algorithm determines the collection of all irreducible compatible team combinations, denoted as $\mathcal{A}_i$, for each task $i$ belonging to the set $\mathcal{N}^{(t)}_{or}$.

\subsubsection{Random task insertion (RaI)}

The algorithm begins by sorting the tasks in $\mathcal{N}^{(t)}_{or}$ in non-increasing order of their priority. Starting from the first task, denoted as $i^*$, in the sorted list, a random team combination is selected from the set $\mathcal{A}_{i^{*}}$. Subsequently, the task is inserted into the optimal positions within the teams comprising the selected random team combination. In the event that $\mathcal{A}_{i^{*}}$ is empty, the task is outsourced.

\subsubsection{Greedy task insertion (GI)}

Similar to random insertion, the algorithm begins by sorting the tasks in $\mathcal{N}^{(t)}_{or}$ in non-increasing order based on their priority. Subsequently, starting from the first task, denoted as $i^*$, in the sorted list, the algorithm selects the optimal team combination from the set $\mathcal{A}_{i^{*}}$, which yields the most favorable insertion cost. The task is then strategically inserted into the optimal positions within the teams belonging to the selected team combination. In the event that $\mathcal{A}_{i^{*}}$ is empty, the task is designated for outsourcing.

\subsubsection{$n$-regret insertion (ReI-$n$)}

Regret heuristics are kind of greedy heuristics enhanced with a look-ahead feature \citep{PotvinRousseau1993}. Based on \citet{PisingerRopke2007}, the regret value for each task $i\in \mathcal{N}^{(t)}_{or}$ is calculated as follows. Let $\Delta \text{TTWT}_k(i)$ denote the change in the objective function value for inserting $i$ at its best team combination instead of its $k^{\textrm{th}}$-best team combination. In the basic version where $n=2$, referred to as \textit{ReI-2}, the task to be inserted next is selected where the cost difference between inserting it into its best team combination and its second best team combination is largest. Then, in each iteration,  depending on the regret heuristic used, i.e., the value of $n$, the regret heuristic chooses the next task $i^*$ to be inserted according to Equation~\ref{Regreti*}:
\begin{align}
i^* = \arg\max_{i\in \mathcal{N}^{(t)}_{or}} \Bigg\{ \sum_{k=2}^{\text{min}\{n,|\mathcal{A}_i|\}}\big(\Delta \text{TTWT}_k(i) - \Delta \text{TTWT}_1(i)\big) \Bigg\} \label{Regreti*}
\end{align}

\subsection{Acceptance criteria }

The acceptance criteria of the simulated annealing is employed in the proposed ALNS. It implies that an improving solution $x^{\text{candidate}}$, where $\text{TTWT}(x^{\text{candidate}})<\text{TTWT}(x^{\text{current}})$, is  accepted with probability 1. If $\text{TTWT}(x^{\text{candidate}}) \geq \text{TTWT}(x^{\text{current}})$, then  $x^{\text{candidate}}$ is accepted with probability  $e^{\frac{-\Delta \text{TTWT}}{T_{iter}}}$, where $\Delta \text{TTWT} = \text{TTWT}(x^{\text{candidate}}) - \text{TTWT}(x^{\text{current}})$ and $T_{iter}$ denotes the temperature of the current iteration. 

\vspace{2mm}

The algorithm initiates with an initial temperature value represented as $T_{0}$. After accepting a non-improving solution, the temperature value is updated by multiplying it by a factor $\alpha$, denoted as $T_{iter} = \alpha T_{iter}$. To reduce the likelihood of accepting non-improving solutions as the algorithm progresses, the temperature value must gradually decrease. Consequently, the parameter $\alpha$ should fall within the range of 0 to 1, i.e., $\alpha \in (0,1)$.

\subsection{Initialization and update of weights}

The initial weights of the heuristics in $\mathcal{D}$ and $\mathcal{R}$ are set to $1/|\mathcal{D}|$ and $1/|\mathcal{R}|$, respectively. During each iteration, the weights of the chosen destroy and repair heuristics, denoted as $h^{\text{destroy}}$ and $h^{\text{repair}}$, respectively, are updated based on the solution quality of the resulting solution $x^{\text{candidate}}$. If $x^{\text{candidate}}$ outperforms the best-found solution, replacing $x^{\text{best}}$, the weights of $h^{\text{destroy}}$ and $h^{\text{repair}}$ are increased by $\sigma_1$. Conversely, if $x^{\text{candidate}}$ improves upon the current solution, replacing $x^{\text{current}}$, the weights of $h^{\text{destroy}}$ and $h^{\text{repair}}$ are increased by $\sigma_2$. However, if $x^{\text{candidate}}$ is accepted based on the acceptance criteria, even without providing a better solution, the weights of $h^{\text{destroy}}$ and $h^{\text{repair}}$ are increased by $\sigma_3$. In cases where $x^{\text{candidate}}$ does not replace $x^{\text{current}}$, the weights of $h^{\text{destroy}}$ and $h^{\text{repair}}$ are increased by $\sigma_4$. It is important to note that the values of $\sigma_1, \sigma_2, \sigma_3$, and $\sigma_4$ are hyperparameters that need to be tuned through computational experiments. However, they must adhere to the following relation:

\begin{equation*}
\sigma_1 > \sigma_2 > \sigma_3 > 0 > \sigma_4
\end{equation*}

At the end of each iteration, the weights of all repair and destroy heuristics are normalized. To choose a destroy and a repair heuristic from the sets $\mathcal{D}$ and $\mathcal{R}$, respectively, we employ the \textit{roulette wheel selection} approach. This approach allocates each alternative an angle proportional to its weight in the roulette wheel.

\subsection{Stopping condition }
The maximum number of iterations, referred to as $\nu_1$, and maximum number of non-improving iterations, referred to as $\nu_2$ are both used as the stopping condition.

\section{Computational experiments}
\label{section:compStudy}

In this section, we evaluate the effectiveness of our proposed reoptimization framework through a series of computational experiments conducted on randomly generated DWSRP instances. The process for creating these instances is detailed in Section \ref{subsec:data}. Our evaluation begins with a comprehensive comparison involving the proposed ALNS algorithm, the state-of-the-art solver CPLEX, and a construction heuristic designed to mimic the decision-making approach of a prudent decision-maker. This comparison is carried out using problem instances representing single periods, as described in Section \ref{section:compStudy_HeuristicMMComparison}. Subsequently, in Section \ref{section:compStudy_largeSizeDWSRPinstances}, we extend our examination to encompass larger instances that span multiple reoptimization periods. This section also explores the impact of various factors, such as the degree of dynamism, effective degree of dynamism, and the number of tasks requiring synchronization, on TWTT. Additionally, we conduct an analysis of the influence of different reoptimization strategies and the effects of varying frozen period durations on solution quality. These aspects are comprehensively explored in Sections~\ref{section:compStudy_reoptimizationPeriodComparison} and ~\ref{section:compStudy_frozenPeriodComparison}, respectively.

\vspace{2mm}

We employed CPLEX 20.1 to solve the proposed mathematical model, setting a time limit of 15 minutes for each run. The implementation of the proposed ALNS-based reoptimization framework was developed in Java. The computational experiments were conducted on a computer equipped with an Intel® Core™ i7-6500U 2.5 GHz processor, complemented by 12 GB of RAM, and running on the Windows 10 operating system. To optimize the performance of the ALNS algorithm, we engaged in a parameter tuning process. As a result, the ALNS algorithm's parameters were meticulously configured to align with the specifications detailed in Table~\ref{tab:ALNSparams}.

\begin{table}[!ht]
\centering
\small
\begin{tabular}{lc}
  \toprule  
 Parameter & Value \\
    \midrule
initial temperature value, $T_0$ & 1000  \\
temperature update coefficient, $\alpha$ & 0.95  \\
max. no. of iterations, $\nu_1$ & 250  \\
max. no. of non-improving iterations, $\nu_2$ & 50 \\
min. degree of destruction, $\gamma^{min}$ &  0.50\\
max. degree of destruction, $\gamma^{max}$ &  1.00\\
weight increase for the best-found improving solution, $\sigma_1$&  0.08\\
weight increase for the current improving solution, $\sigma_2$& 0.05\\
weight increase for the non-improving accepted solution, $\sigma_3$& 0.01\\
weight increase for the non-improving unaccepted solution, $\sigma_4$& -0.03\\
    \bottomrule
\end{tabular}
\caption{Parameter values of the proposed ALNS}
\label{tab:ALNSparams}
\end{table}

\subsection{Data}\label{subsec:data}

In our study, a DWSRP instance represents the problem to be optimized during reoptimization. To distinguish the instance that spans an entire working day, we refer to it as a DWSRP super-instance. An arbitrary DWSRP super-instance, designated as $I(N,K,\delta)$, encompasses a total of $N$ tasks scheduled throughout the entire working day, involving $K$ teams, and embodies a specific degree of dynamism quantified by $\delta$. Our aim is to create instances that manifest diverse degrees of dynamism, precisely $\delta = 0.2, 0.4, 0.6, 0.8$, while varying the total task counts and team numbers, represented as $(N, K) = (30, 2), (40, 3), (50, 4), (60, 5), (75, 5)$.

The process of generating an instance denoted as $I(N,K,\delta)$ unfolds as follows. Considering a standard working day spanning 9 hours, equivalent to 540 minutes, we set $\tau^{\max}=540$. To ensure a gradual distribution of new task arrivals across the work period, we divide the day, denoted as $(0, \tau^{\max}]$, into $d$ time intervals. At the start of the day, $\lfloor N/d\rfloor$ tasks are assumed to be predetermined, thereby establishing the degree of dynamism, $\delta$, for the given instance as $\delta = \big(N-\lfloor N/d\rfloor\big)/N$. The arrival time $a_i$ and the earliest start time $e_i$ for these static tasks are initialized to 0. For each interval $t=1,..., d-1$, corresponding to the time interval $\big(\frac{\tau^{\max}}{d}(t-1), \frac{\tau^{\max}}{d}t\big]$, $\lfloor N/d\rfloor$ tasks are generated. The remaining $\big(N-\lfloor N/d\rfloor (d-1)\big)$ tasks are allocated to the final time interval, i.e., time interval $t=d$, corresponding to $\big(\frac{\tau^{\max}}{d}(d-1), \tau^{\max}]$. In the case of a task $i$ generated for a time interval $t$, its arrival time $a_i$ is determined as a random value within the boundaries of the corresponding interval, while the earliest start time $e_i$ is set to match $a_i$. The latest start time $l_i$ for task $i$ is calculated as the minimum of either $a_i + \textit{Uniform}(10, 50)$ or $\tau^{\max}$.  Furthermore, both the processing time $p_i$ and the priority $w_i$ are established using uniform distributions of $\textit{Uniform(5,25)}$ and $\textit{Uniform(1,5)}$, respectively.
Additionally, each task is assigned a random location within a rectangular metric space measuring 25 km by 25 km. We assume a traveling speed of 30 km per hour for each team, utilizing rectilinear distances to generate the distance matrix.
\vspace{2mm}

\vspace{2mm}

We consider a skill set comprising five distinct skills, denoted as $|\mathcal{Q}|=5$. To determine the skill requirements for each task-skill pair $(i,q)$, we employ a random number generator that generates values between 0 and 1. Specifically, we set $u_{iq}$ to 1 if the randomly generated number falls within the range of 0 to 0.5. This means that task $i$ requires skill $q$ when the random number is less than or equal to 0.5. In cases where a task $i$ does not initially require any skill, i.e., $u_{iq}=0$ for all $q\in \mathcal{Q}$, we then randomly select a skill to be associated with the task. A similar process is applied to determine the skill capabilities of the teams, resulting in the $[v_{kq}]$ matrix. It's worth noting that, depending on the skill requirements of tasks that are not known in advance and the skill capabilities of the teams, some tasks within each super-instance require synchronization. The number of tasks necessitating synchronization is denoted by $N^{\text{sync}}$.

\vspace{2mm}

In this manner, we have generated twenty DWSRP super-instances, each characterized by varying levels of dynamism. Table~\ref{tab:InstanceInfo} provides information on the effective degree of dynamism, $\delta^e$, effective degree of dynamism with time windows, $\delta^e_{tw}$, and the number of tasks requiring synchronization, $N^{\text{sync}}$, for each instance.

\begin{table}[!ht]
\centering
\small
\begin{tabular}{cccc|cccc}
  \toprule  
\textbf{DWSRP-TW-SC}  & \multirow{3}{*}{$\delta^e$} & \multirow{3}{*}{$\delta^e_{tw}$} & \multirow{3}{*}{$N^{\text{sync}}$} & \textbf{DWSRP-TW-SC}  & \multirow{3}{*}{$\delta^e$} & \multirow{3}{*}{$\delta^e_{tw}$} & \multirow{3}{*}{$N^{\text{sync}}$}  \\
\textbf{super-instance}  & & & & \textbf{super-instance}  & &   \\
\textbf{ $I(N,K,\delta)$}  &  & &  & \textbf{ $I(N,K,\delta)$}  &  & \\
\hline
$I(30,2, 0.2)$ & 0.47 & 0.76 & 12& $I(50,4, 0.6)$ & 0.49 & 0.72 &13 \\
$I(30,2, 0.4)$ & 0.77 & 0.72 & 17& $I(50,4, 0.8)$ & 0.53 & 0.7 &19 \\
$I(30,2, 0.6)$ & 0.36 & 0.69 & 13& $I(60,5, 0.2)$ & 0.67 & 0.72 &21 \\
$I(30,2, 0.8)$ & 0.47 & 0.7 & 9& $I(60,5, 0.4)$ &  0.34& 0.68 &11 \\
$I(40,3, 0.2)$ & 0.23 & 0.69 & 17& $I(60, 5, 0.6)$ & 0.91&0.73 &18 \\
$I(40,3, 0.4)$ & 0.25 & 0.72 & 13& $I(60, 5, 0.8)$ & 0.44&0.7 & 12 \\
$I(40,3, 0.6)$ & 0.47 & 0.69& 26& $I(75,5, 0.2)$ & 0.61 &0.74 & 27 \\
$I(40,3, 0.8)$ & 0.62 & 0.71& 8& $I(75,5, 0.4)$ &  0.45 &0.7 & 19  \\
$I(50,4, 0.2)$ & 0.82 & 0.7& 24& $I(75, 5, 0.6)$ & 0.39 &0.71 & 31 \\
$I(50,4, 0.4)$ & 0.12 & 0.7& 7& $I(75, 5, 0.8)$ & 0.18 &0.72 & 12  \\
  \bottomrule
\end{tabular}
\caption{Information about the generated DWSRP super-instances}
\label{tab:InstanceInfo}
\end{table}

\vspace{2mm}

Furthermore, to comprehensively assess the performance of the proposed ALNS algorithm in comparison to both the Mathematical Model (MM) and the Constructive Heuristic (CH), we meticulously crafted a set of thirty small- to medium-size DWSRP instances. These instances encompass a variety of task counts ($4 \leq |\mathcal{N}| \leq 8$), team quantities ($2 \leq |\mathcal{K}| \leq 5$), and skill numbers ($3 \leq |\mathcal{Q}| \leq 5$). The remaining parameters for these instances were generated following the aforementioned methodology. This collection of instances is instrumental in the  analysis presented in Section~\ref{section:compStudy_HeuristicMMComparison}. The complete dataset is available at \citet{mendeley_data}.

\subsection{Analysis on the performance of the proposed ALNS compared to the mathematical model and the constructive heuristic} \label{section:compStudy_HeuristicMMComparison}

To comprehensively investigate the effectiveness of our proposed ALNS approach, we conducted a thorough comparative analysis, involving the ALNS, the Mathematical Model (abbreviated as MM), and the Constructive Heuristic (abbreviated as CH). This evaluation involved multiple reoptimization problems across the range of small- to medium-sized instances. Comprehensive outcomes of these comparative experiments can be found in Table~\ref{tab:compStudy_HeuristicMMComparison}. This table provides a comprehensive overview of each DWSRP instance, including essential parameters such as task count ($|\mathcal{N}|$), team count ($|\mathcal{K}|$), and skill count ($|\mathcal{Q}|$). It further presents the objective function value, denoted as TWTT, as computed by the MM, CH, and ALNS techniques. Furthermore, the table highlights the percentage gap associated with the best found solution by the MIP within limited time (15 minutes) (under column \textit{MIP Gap (\%)}), the percentage of improvement in the objective value provided by the ALNS compared to the MIP (under column \textit{ALNS Diff (\%)}), the percentage of improvement in the objective value provided by the ALNS compared to CH (under column \textit{ALNS Imp (\%)}), and the respective runtimes of MM, CH, and ALNS.
\begin{landscape}
\begin{table}[htbp]
\centering
\small
\begin{tabular}{c c c c cccccc c ccc}
  \toprule  
\textbf{ DWSRP-TW-SC}  & \multirow{3}{*}{$N$}  & \multirow{3}{*}{$K$} & \multirow{3}{*}{$Q$} &
\multicolumn{6}{c}{TWTT} & & \multicolumn{3}{c}{CPU(s)}  \\
\cline{5-10} \cline{12-14}
\multirow{2}{*}{\textbf{ Instance}} &   &  &  & \multirow{2}{*}{MM} & \multirow{2}{*}{CH} & \multirow{2}{*}{ALNS} & MIP & ALNS & ALNS  &&  \multirow{2}{*}{MM} & \multirow{2}{*}{CH} & \multirow{2}{*}{ALNS} \\
  &   &  &  &  &  &  &  Gap (\%) &  Diff (\%) &  Imp (\%) &&   &  & 
  \\
  \hline
1&4&2&3&769&769&769&0.00&0.00&0.00&&$<1$&$<1$&$<1$ \\
2&4&3&3&747&747&747&0.00&0.00&0.00&&$<1$&$<1$&$<1$ \\
3&4&4&4&608&608&608&0.00&0.00&0.00&&$<1$&$<1$&1\\
4&5&2&3&1335&1335&1335&0.00&0.00&0.00&&$<1$&$<1$&1 \\
5&5&3&3&1423&1467&1423&0.00&0.00&3.00&&$<1$&$<1$&2 \\
6&5&4&4&1141&1165&1141&0.00&0.00&2.06&&$<1$&$<1$&1 \\
7&5&5&5&909&1009&909&0.00&0.00&9.91&&$<1$&$<1$&4 \\
8&6&2&3&1697&1716&1697&0.00&0.00&1.11&&$<1$&$<1$&1 \\
9&6&3&4&1702&2906&1702&0.00&0.00&41.43&&$<1$&$<1$&3 \\
10&6&5&3&1222&1276&1222&0.00&0.00&4.23&&5&$<1$&3 \\
11&7&2&3&4601&4601&4601&0.00&0.00&0.00&&$<1$&$<1$&$<1$ \\
12&7&3&3&5232&5232&5232&0.00&0.00&0.00&&$<1$&$<1$&$<1$\\
13&7&4&4&1669&2923&1669&0.00&0.00&42.90&&$<1$&$<1$&4 \\ 
14&7&5&5&1556&2290&1556&0.00&0.00&32.05&&4&$<1$&6 \\ 
15&8&2&5&4685&4685&4685&0.00&0.00&0.00&&$<1$&$<1$&4 \\
16&8&3&3&1892&2070&1892&0.00&0.00&8.60&&16&$<1$&3 \\
17&8&4&4&1396&1657&1396&0.00&0.00&15.75&&109&$<1$&5\\ 
18&10&2&3&3685&3817&3685&0.00&0.00&3.46&&207&$<1$&5 \\ 
19&10&3&3&5271&5601&5271&0.00&0.00&5.89&&120&$<1$&10 \\
20&10&4&4&2155&2946&2120&53.51&1.62&28.04&&$ 900$&2&6\\
21&10&5&5&17621&17621&17621&0.00&0.00&0.00&&$<1$&$<1$&$<1$\\
22&12&2&3&10779&10782&10779&0.00&0.00&0.03&&$<1$&$<1$&2 \\ 
23&12&3&4&4473&4647&4400&71.03&1.63&5.32&&$900$&$<1$&9 \\ 
24&12&5&5&3382&3473&3105&69.50&8.19&10.60&&$ 900$&8&20\\
25&15&2&3&9275&9985&9275&0.00&0.00&7.11&&72&$<1$&5\\ 
26&15&3&4&9314&9874&8561&89.75&8.08&13.30&&$ 900$&$<1$&27 \\ 
27&15&5&5&5890&5721&3783&84.77&35.77&33.88&&$ 900$&10&47 \\
28&20&2&3&21376&21376&21376&0.00&0.00&0.00&&$<1$&$<1$&3 \\
29&20&3&4&5991&5925&5925&82.22&1.10&0.00&&$ 900$&5&29\\  
30&20&5&5&6797&5212&4908&87.05&27.79&5.83&&$900$&11&85\\
  \bottomrule
  &   &  &  &  &  & \textbf{Min.} &  0.00 &  35.77 &  0.00 && 1  &  & 1
  \\
    &   &  &  &  &  & \textbf{Max.} &  89.75 &  0.00 &  42.90 && 900  &  & 85
  \\
    &   &  &  &  &  & \textbf{Avg.} &  17.93 &  2.81 &  9.15 && 228.3  &  & 9.7
  \\
  \bottomrule
\end{tabular}
\caption{Comparing ALNS, MM, and CH Results for small-to medium-size DWSRP instances }
\label{tab:compStudy_HeuristicMMComparison}
\end{table}

\end{landscape}
\vspace{2mm}

As presented in Table~\ref{tab:compStudy_HeuristicMMComparison}, among the 30 problem instances, the MM was able to attain optimal solutions (as indicated by zero percentage gap values in  the MIP Gap ($\%$)  column) for 23 instances within the prescribed 15-minute runtime limit. Remarkably, for these instances, our propsoed ALNS also achieved the optimal solution, often within the same or significantly less time. For the remaining 7 instances where the MM couldn't reach optimality withing the limited runtime, a comparison of MM and ALNS objective function values (found under the ALNS Diff ($\%$)  column) reveals ALNS consistently outperforming the MM.   Notably, the ALNS demonstrated an average improvement of $12.03\%$ over the MM.  A closer look at the average runtime reinforces ALNS's superiority, as it proved to be approximately 23.5 times faster than the MM.  Taking both the quality of objective function results and the efficiency of the ALNS algorithm into account, it's evident that ALNS significantly outshines the MM.

\vspace{2mm}

Turning to the values in the ALNS Imp ($\%$) column, we further assess the extent to which the proposed ALNS enhances solutions obtained from the CH. This evaluation mirrors ALNS's practical value, given the preference for simple greedy algorithms in many industries. Upon analyzing the results, a consistent trend emerges: ALNS consistently improves the objective function value, reducing TWTT values by an average of $9.15\%$. Notably, some instances saw substantial improvements, reaching as high as $42.90\%$. It's worth highlighting that ALNS achieves these enhancements in an average time of just 9.7 seconds.

\subsection{Analysis on the results of the proposed ALNS  on the large-size DWSRP super-instances } \label{section:compStudy_largeSizeDWSRPinstances}

Next, we apply both the CH and the ALNS independently within the proposed reoptimization framework for each large-size DWSRP super-instance. Following an initial analysis, we set the values of $\beta^{\text{task}}$ and  $\beta^{\text{time}}$ that indicate the frequency of the reoptimizations at 5 tasks and 60 time units, respectively, while the frozen period length is defined as $f=30$ time units.  For each approach employed on a DWSRP super-instance, we aggregate the results of corresponding DWSRP instances across each reoptimization cycle. Subsequently, we  compute the total waiting time of all tasks within the given super-instance, presenting this value under the relevant $\text{TWTT}$ column in Table~\ref{tab:compStudy_DWSPsuperinstances}.

\vspace{2mm}

The number of reoptimizations performed during the workday remains consistent for both the CH and ALNS approaches across each DWSRP super-instance, and this count is documented in the \textit{No. of reopt.} column. Additionally, for each DWSRP super-instance, we present the improvement percentage attributed to the use of ALNS compared to CH within the reoptimization framework. Furthermore, we include the total runtimes of CH and ALNS, taking into account reoptimizations conducted for the corresponding instances.

\begin{table}[htbp]
\centering
\small
\begin{tabular}{c ccc c cc}
  \toprule  
\textbf{ DWSRP-TW-SC}  & 
\multicolumn{3}{c}{TWTT} & \multirow{3}{*}{\textbf{No. of reopt.}}  & \multicolumn{2}{c}{CPU(s)}  \\
\cline{2-4} \cline{6-7}
\multirow{2}{*}{\textbf{Super-instance}} &   \multirow{2}{*}{CH} & \multirow{2}{*}{ALNS} & ALNS & &  \multirow{2}{*}{CH} & \multirow{2}{*}{ALNS} \\
  &   &  &   Imp (\%) & &   &  
  \\
  \hline
$I(30,2,0.2)$&18290&14445&21.02&8&0.57&42.12\\
$I(30,2,0.4)$&15843&15135&4.47&13&0.31&24.62\\
$I(30,2,0.6)$&9496&7837&17.47&19&0.68&31.49\\
$I(30,2,0.8)$&7860&7558&3.84&25&0.29&10.66\\
$I(40,3,0.2)$&19857&16504&16.89&9&7.53&308.58\\
$I(40,3,0.4)$&13387&13110&2.07&17&7.81&511.49\\
$I(40,3,0.6)$&14015&13139&6.25&24&3.44&276.91\\
$I(40,3,0.8)$&7378&6792&7.94&32&0.94&15.76\\
$I(50,4,0.2)$&26199&23595&9.94&11&161.71&453.01\\
$I(50,4,0.4)$&14658&11818&19.38&21&41.99&619.32\\
$I(50,4,0.6)$&17466&11901&31.86&31&34.86&685.3\\
$I(50,4,0.8)$&10271&9752&5.05&41&9.48&202.3\\
$I(60,5,0.2)$&26991&24638&8.72&13&716.36&1386.44\\
$I(60,5,0.4)$&25604&20290&20.75&25&217.76&3078.96\\
$I(60,5,0.6)$&18151&14988&17.43&37&141.94&855.13\\
$I(60,5,0.8)$&17672&11237&36.41&49&26.76&263.43\\
$I(75,5,0.2)$&44429&36635&17.54&16&3526.19&31996.16\\
$I(75,5,0.4)$&27496&23380&14.97&31&678.66&1716.18\\
$I(75,5,0.6)$&28462&21757&23.56&45&137.53&1204.65\\
$I(75,5,0.8)$&26611&17598&33.87&59&71.69&370.38\\
  \bottomrule
\textbf{Min.} &7378&6792&2.07&8&0.29&10.66\\

\textbf{Max.} &44429&36635&36.41&59&3526.19&31996.16\\

\textbf{Avg.} &19506.8&16105.45&15.97&26.3&289.33&2202.64\\

  \bottomrule
\end{tabular}
\caption{The results of the reoptimization frameworks that employ the CH and the ALNS, separately,  on large-size DWSRP super-instances }
\label{tab:compStudy_DWSPsuperinstances}
\end{table}

\vspace{2mm}

The data presented in Table~\ref{tab:compStudy_DWSPsuperinstances} reveals a substantial average reduction of nearly $16\%$ in TWTT, with the potential for this improvement to extend to almost $37\%$ in specific instances. However, it is evident from these findings that achieving such a significant substantial enhancement comes at the cost of increased solution times. This results in ALNS requiring approximately $7.60$ times the CPU time compared to the construction heuristic. These results  suggest that the viability of adopting ALNS depends on the nature of the business models; decision-makers must carefully evaluate the trade-off between the computational time consumed by ALNS and its potential benefits. Nonetheless, the results underscore that if the algorithm's time requirements are accommodated, decision-makers can gain a significant advantage from implementing ALNS.

\vspace{2mm}

\begin{figure}[ht!]
  \centering
  \includegraphics[width=0.6\linewidth]{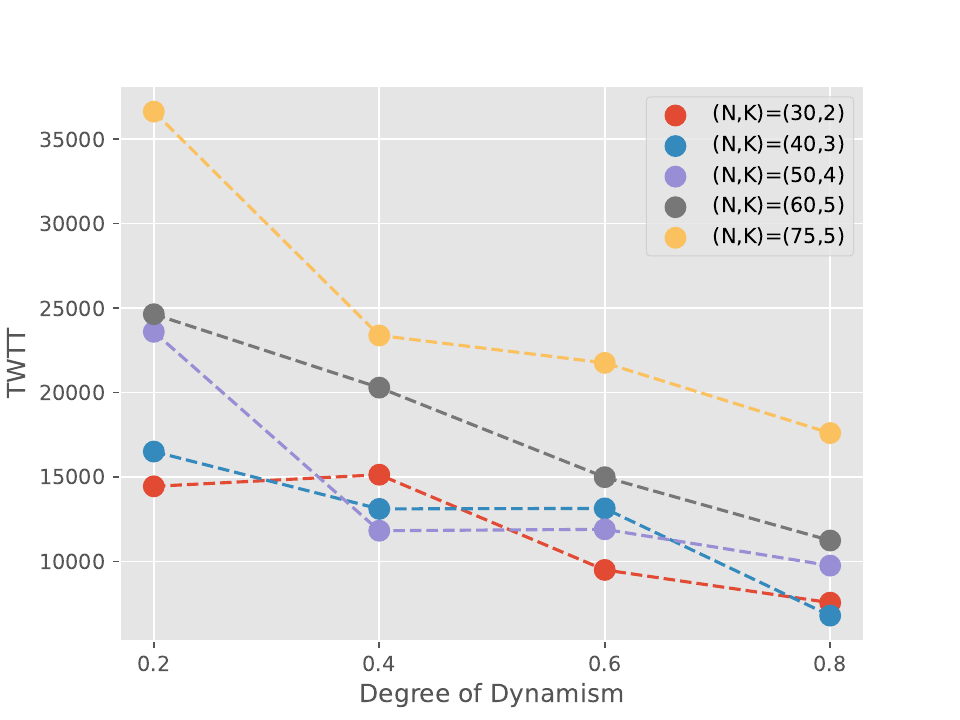}
  \caption{Impact of degree of dynamism ($\delta$) on TWTT}
  \label{fig:delta_twtt}
\end{figure}

\begin{figure}[ht!]
  \centering
  \includegraphics[width=1.1\linewidth]{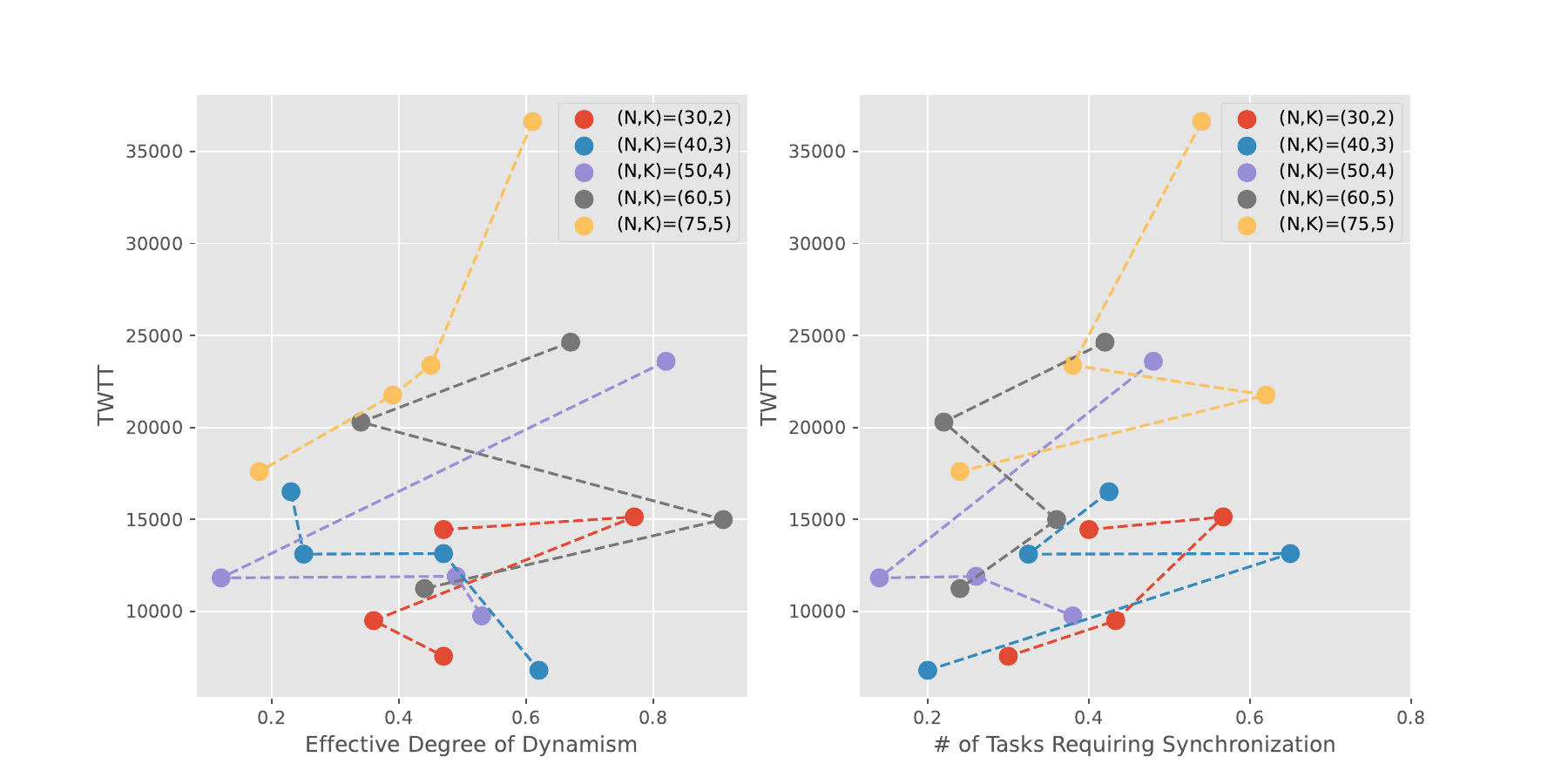}
  \caption{\textbf{(a)} Impact of effective degree of dynamism ($\delta^e$) on TWTT \textbf{(b)} Impact of $N^{sync}$ on TWTT}
  \label{fig:useless}
\end{figure}

Table~\ref{tab:compStudy_DWSPsuperinstances} offers valuable insights into the influence of several factors associated with super-instances. These factors encompass the level of dynamism, effective level of dynamism, and the number of tasks requiring synchronization, all of which have the potential to affect TWTT. As depicted in Figure~\ref{fig:delta_twtt}, it is evident that TWTT tends to decrease with an increase in the level of dynamism ($\delta$). This observation aligns with our expectations, as a higher influx of tasks during execution diminishes the flexibility of the planning framework. However, as illustrated in Figure~\ref{fig:useless}, there does not exist a direct correlation between $\delta^e$-TWTT and $N^{sync}$-TWTT. This lack of correlation may arise from neither $\delta^e$ nor $N^{sync}$ individually having a significant impact on TWTT. Since all our datasets exhibit similar levels of $\delta^e_{tw}$ values, and it's worth noting that $\delta^e_{tw}$ is essentially an extension of $\delta^e$, we chose not to delve into the specific impact of $\delta^e_{tw}$ in this analysis.

\subsection{Analysis on different reoptimization strategies} \label{section:compStudy_reoptimizationPeriodComparison}

The aim of this section is to investigate the impact of employing various reoptimization strategies on both the objective function value (TWTT) and the total CPU time. These strategies are defined by different values of $\beta^{\text{task}}$. To achieve this objective, we analyze two distinct sets of super-instances:
(i) the super-instances listed in Table \ref{tab:InstanceInfo}, referred to as \textit{tight instances}, and
(ii) super-instances identical to the tight instances, with the exception of task deadlines. This second set, referred to as \textit{loose instances}, uniformly assigns task deadlines at $\tau^{\max}$. The inclusion of these additional super-instances eliminates potential biases in the results that could arise from  tasks with particularly stringent deadlines necessitating outsourcing.  

\vspace{2mm}

In this experiment, we initially focus on the tight super-instances. For each individual super-instance within this tight set, we conduct separate simulations of the entire workday, utilizing the options $\beta^{\text{task}}=1$, $\beta^{\text{task}}=3$, and $\beta^{\text{task}}=5$ each with a frozen period length of zero. Our ALNS algorithm is employed for each reoptimization, and we present the resulting TWTT values, CPU times, and the cumulative count of outsourced tasks in Table \ref{tab:compStudy_tight_reoptimization}. Subsequently, a parallel procedure is  executed for the loose super-instances, with the corresponding outcomes documented in Table \ref{tab:compStudy_looose_reoptimization}.

\vspace{2mm}

\begin{table}[htbp]
\centering
\scriptsize
\begin{tabular}{c ccc c c ccc c c ccc}
  \toprule  
\textbf{Tight}  & 
\multicolumn{3}{c}{$\beta^{\text{task}}=1$}   & 
\multicolumn{1}{c}{}   & 
\multicolumn{3}{c}{$\beta^{\text{task}}=3$}   &
\multicolumn{1}{c}{}   & 
\multicolumn{3}{c}{$\beta^{\text{task}}=5$}  \\
\cline{2-4} \cline{6-8}  \cline{10-12}
\multirow{2}{*}{\textbf{Super-instance}} &   \multirow{2}{*}{TWTT} & \multirow{2}{*}{CPU(s)} & \multirow{2}{*}{outso.} & & \multirow{2}{*}{TWTT} & \multirow{2}{*}{CPU(s)} & \multirow{2}{*}{outso.} & &\multirow{2}{*}{TWTT} & \multirow{2}{*}{CPU(s)} & \multirow{2}{*}{outso.} \\
  \\
  \hline
$I(30,2,0.2)$&14445&42.12&14&&15369&21.03&15&&17001&18.31&14\\
$I(30,2,0.4)$&15135&24.62&13&&17398&19.11&16&&18087&18.66&12\\
$I(30,2,0.6)$&7837&31.49&7&&8872&10.48&8&&12389&13.02&5\\
$I(30,2,0.8)$&7558&10.66&8&&9893&3.3&8&&12764&5.33&4\\
$I(40,3,0.2)$&16504&308.58&10&&20972&199.3&11&&22589&411.14&11\\
$I(40,3,0.4)$&13110&511.49&9&&15347&360.39&10&&18922&366.77&14\\
$I(40,3,0.6)$&13139&276.91&8&&13548&157.98&8&&18996&172.39&6\\
$I(40,3,0.8)$&6792&15.76&5&&9710&5.86&5&&11710&11.09&6\\
$I(50,4,0.2)$&23595&453.01&14&&27781&787.1&13&&27840&688.79&16\\
$I(50,4,0.4)$&11818&619.32&9&&14893&129.08&8&&15377&220.63&13\\
$I(50,4,0.6)$&11901&685.3&10&&12705&336.97&10&&14856&199.87&14\\
$I(50,4,0.8)$&9752&202.3&7&&10663&42.46&7&&13037&32.72&10\\
$I(60,5,0.2)$&24638&1386.44&17&&25404&1656.26&21&&27959&896.53&18\\
$I(60,5,0.4)$&20290&3078.96&12&&24120&846.34&16&&25418&664.03&11\\
$I(60,5,0.6)$&14988&855.13&9&&15367&290.88&11&&16205&250.99&9\\
$I(60,5,0.8)$&11237&263.43&9&&12233&131.91&7&&14303&124.38&5\\
$I(75,5,0.2)$&36635&31996.16&29&&39205&16444.81&33&&42754&13702.6&26\\
$I(75,5,0.4)$&23380&1716.18&19&&23630&1305.72&22&&24988&1014.62&21\\
$I(75,5,0.6)$&21757&1204.65&22&&23347&862.76&19&&25579&510.55&27\\
$I(75,5,0.8)$&17598&370.38&12&&18827&179.95&13&&24423&120.7&11\\
  \bottomrule
  \textbf{Min.} &6792&10.66&5&&8872&3.3&5&&11710&5.33&4\\
  \textbf{Max.}  &36635&31996.16&29&&39205&16444.81&33&&42754&13702.6&27\\
  \textbf{Avg.} &16105.45&2202.64&12.15&&17964.2&1189.58&13.05&&20259.85&972.16&12.65\\
  \bottomrule
\end{tabular}
\caption{Results on tight super-instances with changing values of $\beta^{\text{task}}$}
\label{tab:compStudy_tight_reoptimization}
\end{table}

\begin{table}[htbp]
\centering
\scriptsize
\begin{tabular}{c ccc c c ccc c c ccc}
  \toprule  
\textbf{Loose}  & 
\multicolumn{3}{c}{$\beta^{\text{task}}=1$}   & 
\multicolumn{1}{c}{}   & 
\multicolumn{3}{c}{$\beta^{\text{task}}=3$}   &
\multicolumn{1}{c}{}   & 
\multicolumn{3}{c}{$\beta^{\text{task}}=5$}  \\
\cline{2-4} \cline{6-8}  \cline{10-12}
\multirow{2}{*}{\textbf{Super-instance}} &   \multirow{2}{*}{TWTT} & \multirow{2}{*}{CPU(s)} & \multirow{2}{*}{outso.} & & \multirow{2}{*}{TWTT} & \multirow{2}{*}{CPU(s)} & \multirow{2}{*}{outso.} & &\multirow{2}{*}{TWTT} & \multirow{2}{*}{CPU(s)} & \multirow{2}{*}{outso.} \\
  \\
  \hline
$I(30,2,0.2)$&12845&81.59&9&&13089&52.19&9&&13991&126.33&9\\
$I(30,2,0.4)$&13854&44.44&10&&14904&16.47&11&&16147&21.94&8\\
$I(30,2,0.6)$&7837&23.41&7&&9901&11.14&7&&12094&20.79&6\\
$I(30,2,0.8)$&7922&17.86&7&&8088&7.49&6&&12473&6.35&4\\
$I(40,3,0.2)$&14747&643.68&2&&16721&459.8&1&&18384&590.16&0\\
$I(40,3,0.4)$&11170&957.53&2&&12357&697.99&4&&13406&385.3&5\\
$I(40,3,0.6)$&13839&432.43&9&&12247&256.95&6&&19372&220.09&4\\
$I(40,3,0.8)$&6552&14.75&4&&9784&10.9&4&&10963&11.91&5\\
$I(50,4,0.2)$&20308&814.25&7&&22113&595.68&7&&23161&1019.99&9\\
$I(50,4,0.4)$&11603&717.35&9&&14249&328.51&8&&12957&183.45&11\\
$I(50,4,0.6)$&10676&826.65&7&&11378&629.85&8&&13050&554.01&12\\
$I(50,4,0.8)$&8624&120.91&6&&10663&36.96&7&&11044&56.17&7\\
$I(60,5,0.2)$&18902&2043.36&5&&19730&1912.13&6&&23213&1809.02&5\\
$I(60,5,0.4)$&18228&8591.31&6&&18362&2211.29&6&&23693&2100.76&2\\
$I(60,5,0.6)$&10780&1005.91&7&&11871&445.09&8&&14785&423.12&9\\
$I(60,5,0.8)$&10318&436.92&8&&11423&99.27&8&&14615&131.27&5\\
$I(75,5,0.2)$&25944&57847.9&6&&26361&35707.8&8&&27275&40546.73&7\\
$I(75,5,0.4)$&20180&4202.78&11&&20342&2131.41&11&&22350&2435.08&12\\
$I(75,5,0.6)$&20513&2327.05&20&&21953&1114.38&16&&20813&985.76&22\\
$I(75,5,0.8)$&17051&604.61&12&&17445&223.78&13&&22037&100.35&8\\
  \bottomrule
  \textbf{Min.} &6552&14.75&2&&8088&7.49&1&&10963&6.35&0\\
  \textbf{Max.}  &25944&57847.9&20&&26361&35707.8&16&&27275&40546.73&22\\
  \textbf{Avg.} &14094.65&4087.73&7.7&&15149.05&2347.45&7.7&&17291.15&2586.43&7.5\\
  \bottomrule
\end{tabular}
\caption{Results on loose super-instances with changing values of $\beta^{\text{task}}$}
\label{tab:compStudy_looose_reoptimization}
\end{table}

Table \ref{tab:compStudy_tight_reoptimization} and Table \ref{tab:compStudy_looose_reoptimization} clearly illustrate that reducing the value of $\beta^{\text{task}}$ leads to improved (lower) TWTT values, regardless of whether the super-instance is classified as tight or loose. This phenomenon arises because a reduced number of new tasks into the system for reoptimization provides greater degrees of freedom to our optimization framework. Consequently, lower values of $\beta^{\text{task}}$ are more likely to yield enhanced TWTT values. This observation remains consistent across all super-instances listed in Table \ref{tab:compStudy_tight_reoptimization}. On average, reducing $\beta^{\text{task}}$ from five to three results in a gain of $11.33\%$, while decreasing it from three to one contributes $10.35\%$ to the improvement. Similarly, this observation holds true in Table \ref{tab:compStudy_looose_reoptimization} for nearly every super-instance, with the exception of one instance. On average across these super-instances, decreasing $\beta^{\text{task}}$ from five to three yields an improvement of $12.38\%$, while reducing it from three to one results in a contribution of $6.96\%$ to the enhancement.  However, this enhancement achieved through a reduced number of required tasks to trigger the subsequent reoptimization phase is accompanied by an increase in the elapsed CPU time. In the context of the tight super-instances, as depicted in Table \ref{tab:compStudy_tight_reoptimization}, diminishing $\beta^{\text{task}}$ from five to three results in an average CPU time increase of $22.36\%,$ while reducing it from three to one elevates the average CPU time by $85.16\%.$ Meanwhile, in the case of the loose instances, the behavior is not evident upon reducing $\beta^{\text{task}}$ from five to three. However, a reduction from three to one amplifies the elapsed CPU time by $74.14\%.$ This substantial increase in solution time is to be expected, given that a lower value of $\beta^{\text{task}}$ corresponds to a higher frequency of invoking the ALNS algorithm.

\vspace{2mm}

It's worth highlighting that the observation we've made—namely, that reducing $\beta^{\text{task}}$ improves the objective function value while increasing CPU time—is consistent for both types of super-instances. Furthermore, no evidence suggests a correlation between the decision of $\beta^{\text{task}}$ and the number of outsourced tasks. Nonetheless, as anticipated, we observe a notable reduction in the number of outsourced tasks for the loose super-instances. This outcome can be attributed primarily to the absence of strict deadlines in these cases. The decision-maker possesses the flexibility to defer task assignment to outsourcing in favor of accommodating them in subsequent reoptimization phases.

\subsection{Analysis on the length of the frozen period, $f$} \label{section:compStudy_frozenPeriodComparison}

Recall that the purpose of implementing a frozen period of suitable duration is to minimize confusion during the transition from one TTP to another. With this goal in mind, this section is dedicated to understanding the influence of varying frozen period lengths, enabling us to offer practical policy recommendations for industry decision-makers. In this context, similar to our approach in Section \ref{section:compStudy_reoptimizationPeriodComparison}, we utilize both tight and loose instances to augment the robustness of our analysis.

\vspace{2mm}

In this experiment, for each individual super-instance within the sets of tight and loose instances, we conduct separate simulations covering an entire workday. We explore three options: setting $f$ to 0 (corresponding to no frozen period), setting $f$ to $\max (t^{\max},p^{\max})$ where $t^{\max}$ and $p^{\max}$ represent the maximum travel time and processing time for the corresponding super-instance, respectively, and setting $f$ to $t^{\max}+p^{\max}$, all with a $\beta^{\text{task}}$ value of 1. It is important to note that  $t^{\max}+p^{\max} > \max (t^{\max},p^{\max}) > 0$, as both $t^{\max}$ and $p^{\max}$ are strictly positive. For reoptimization, we employ our ALNS algorithm and record the resulting values for TWTT, CPU times, and the cumulative count of outsourced tasks. These results are presented in Table \ref{tab:compStudy_frozen_tight} and \ref{tab:compStudy_frozen_loose} for tight and loose super-instances, respectively. 

\begin{table}[htbp]
\centering
\scriptsize
\begin{tabular}{c ccc c c ccc c c ccc}
  \toprule  
\textbf{ Tight}  & 
\multicolumn{3}{c}{f=0}   & 
\multicolumn{1}{c}{}   & 
\multicolumn{3}{c}{$f=\max (t^{\max},p^{\max})$}   &
\multicolumn{1}{c}{}   & 
\multicolumn{3}{c}{$f=t^{\max}+p^{\max}$}  \\
\cline{2-4} \cline{6-8}  \cline{10-12}
\multirow{2}{*}{\textbf{Super-instance}} &   \multirow{2}{*}{TWTT} & \multirow{2}{*}{CPU(s)} & \multirow{2}{*}{outso.} & & \multirow{2}{*}{TWTT} & \multirow{2}{*}{CPU(s)} & \multirow{2}{*}{outso.} & &\multirow{2}{*}{TWTT} & \multirow{2}{*}{CPU(s)} & \multirow{2}{*}{outso.} \\
  \\
  \hline
$I(30,2,0.2)$&14445&42.12&14&&15286&19.75&13&&16132&66.15&14\\
$I(30,2,0.4)$&15135&24.62&13&&15933&31.21&14&&15958&26.66&14\\
$I(30,2,0.6)$&7837&31.49&7&&10779&16.83&9&&11319&9.34&8\\
$I(30,2,0.8)$&7558&10.66&8&&10227&5.36&9&&11314&4.39&9\\
$I(40,3,0.2)$&16504&308.58&10&&18061&319.73&11&&18268&456&10\\
$I(40,3,0.4)$&13110&511.49&9&&14169&276.24&9&&15039&293.08&9\\
$I(40,3,0.6)$&13139&276.91&8&&15502&76.56&9&&15998&55.29&12\\
$I(40,3,0.8)$&6792&15.76&5&&9261&7.61&8&&12767&11.28&10\\
$I(50,4,0.2)$&23595&453.01&14&&24277&514.66&15&&25772&365.09&13\\
$I(50,4,0.4)$&11818&619.32&9&&12693&177.97&9&&13235&230.9&9\\
$I(50,4,0.6)$&11901&685.3&10&&13959&547.54&12&&14620&331.44&10\\
$I(50,4,0.8)$&9752&202.3&7&&15338&41.98&12&&15860&21.99&11\\
$I(60,5,0.2)$&24638&1386.44&17&&23415&1660.26&15&&22736&1028.8&16\\
$I(60,5,0.4)$&20290&3078.96&12&&21427&1345.96&15&&23657&1184.37&13\\
$I(60,5,0.6)$&14988&855.13&9&&16779&374.6&10&&18640&307.44&10\\
$I(60,5,0.8)$&11237&263.43&9&&17391&136.34&13&&18603&57.69&12\\
$I(75,5,0.2)$&36635&31996.16&29&&38376&19140.52&29&&37978&14228.26&29\\
$I(75,5,0.4)$&23380&1716.18&19&&24846&1197.53&19&&24887&976.14&21\\
$I(75,5,0.6)$&21757&1204.65&22&&24637&624.05&23&&25477&587.17&25\\
$I(75,5,0.8)$&17598&370.38&12&&27507&88.07&20&&28748&75.67&21\\
  \bottomrule
  \textbf{Min.} &6792&10.66&5&&9261&5.36&8&&11314&4.39&8\\
  \textbf{Max.}  &36635&31996.16&29&&38376&19140.52&29&&37978&14228.26&29\\
  \textbf{Avg.} &16105.45&2202.64&12.15&&18493.15&1330.14&13.7&&19350.4&1015.86&13.8\\
  \bottomrule
\end{tabular}
\caption{Results on tight super-instances with changing values of $f$}
\label{tab:compStudy_frozen_tight}
\end{table}

\begin{table}[htbp]
\centering
\scriptsize
\begin{tabular}{c ccc c c ccc c c ccc}
  \toprule  
\textbf{ Loose}  & 
\multicolumn{3}{c}{f=0}   & 
\multicolumn{1}{c}{}   & 
\multicolumn{3}{c}{$f=\max (t^{\max},p^{\max})$}   &
\multicolumn{1}{c}{}   & 
\multicolumn{3}{c}{$f=t^{\max}+p^{\max}$}  \\
\cline{2-4} \cline{6-8}  \cline{10-12}
\multirow{2}{*}{\textbf{Super-instance}} &   \multirow{2}{*}{TWTT} & \multirow{2}{*}{CPU(s)} & \multirow{2}{*}{outso.} & & \multirow{2}{*}{TWTT} & \multirow{2}{*}{CPU(s)} & \multirow{2}{*}{outso.} & &\multirow{2}{*}{TWTT} & \multirow{2}{*}{CPU(s)} & \multirow{2}{*}{outso.} \\
  \\
  \hline
$I(30,2,0.2)$&12845&81.59&9&&12836&133.63&9&&12836&52.18&9\\
$I(30,2,0.4)$&13854&44.44&10&&14330&37.2&10&&14422&33.37&10\\
$I(30,2,0.6)$&7837&23.41&7&&9727&21.39&8&&9893&11.94&7\\
$I(30,2,0.8)$&7922&17.86&7&&10229&12.47&7&&10891&6.91&8\\
$I(40,3,0.2)$&14747&643.68&2&&15085&794.35&2&&14879&570.69&2\\
$I(40,3,0.4)$&11170&957.53&2&&12068&701.81&4&&12122&585.81&2\\
$I(40,3,0.6)$&13839&432.43&9&&14426&262.29&9&&15334&139.47&10\\
$I(40,3,0.8)$&6552&14.75&4&&9793&11.86&6&&10863&8.11&8\\
$I(50,4,0.2)$&20308&814.25&7&&20132&944.27&7&&21107&813.64&7\\
$I(50,4,0.4)$&11603&717.35&9&&12242&427.6&9&&12105&450.59&9\\
$I(50,4,0.6)$&10676&826.65&7&&12370&500.91&8&&13307&507.85&7\\
$I(50,4,0.8)$&8624&120.91&6&&14260&32.62&12&&16929&56.3&12\\
$I(60,5,0.2)$&18902&2043.36&5&&19281&2357.89&5&&19971&2379.95&5\\
$I(60,5,0.4)$&18228&8591.31&6&&19530&4555.02&5&&22089&2292.28&7\\
$I(60,5,0.6)$&10780&1005.91&7&&13864&422.13&7&&15544&278.19&7\\
$I(60,5,0.8)$&10318&436.92&8&&16249&141.05&10&&18253&97.08&12\\
$I(75,5,0.2)$&25944&57847.9&6&&24901&52782.11&5&&26071&55740.71&5\\
$I(75,5,0.4)$&20180&4202.78&11&&21920&3106.4&14&&20903&2570.13&10\\
$I(75,5,0.6)$&20513&2327.05&20&&22202&814.07&22&&24499&633.46&24\\
$I(75,5,0.8)$&17051&604.61&12&&23045&102.09&14&&25807&67.96&15\\
  \bottomrule
  \textbf{Min.} &6552&14.75&2&&9727&11.86&2&&9893&6.91&2\\
  \textbf{Max.} &25944&57847.9&20&&24901&52782.11&22&&26071&55740.71&24\\
  \textbf{Avg.} &14094.65&4087.73&7.7&&15924.5&3408.06&8.65&&16891.25&3364.83&8.8\\
  \bottomrule
\end{tabular}
\caption{Results on loose super-instances with changing values of $f$}
\label{tab:compStudy_frozen_loose}
\end{table} 

Tables \ref{tab:compStudy_frozen_tight} and \ref{tab:compStudy_frozen_loose} provide clear evidence that extending the duration of the frozen period has a detrimental effect on TWTT values, regardless of whether the super-instance falls into the tight or loose category. This outcome aligns with expectations, a longer frozen period reduces our algorithm's flexibility, as it narrows the set of tasks the algortihm can optimize. These extended durations, intended to prevent confusion during transitions between different plans, can be seen as a trade-off between flexibility and stability. Specifically, when we extend the frozen period for the tight super-instances from $0$ to $\max (t^{\max},p^{\max})$, we observe an average reduction in TWTT of $14.83\%$. Further increasing it from $\max (t^{\max},p^{\max})$ to $t^{\max}+p^{\max}$ results in an average reduction of $4.64\%$. For loose super-instances, similarly, increasing the frozen period from $0$ to $\max (t^{\max},p^{\max})$ leads to an average TWTT reduction of $12.98\%$, while extending it from $\max (t^{\max},p^{\max})$ to $t^{\max}+p^{\max}$ results in an average reduction of $6.07\%$.

\vspace{2mm}

The  analysis above holds practical relevance for managers considering investments in technologies aimed at reducing the duration of the frozen period. Consider a scenario where costly technology can diminish the value of $f$ from $t^{\max}+p^{\max}$ to $\max (t^{\max},p^{\max})$. In such a case, the adoption of this technology may not be justified, given the expected improvement is only around $5\%$. Conversely, technology capable of reducing the frozen period from $\max (t^{\max},p^{\max})$ to nearly 0 warrants serious consideration, as it promises a significantly substantial enhancement.

\vspace{2mm}

While reducing the duration of the frozen period offers the advantage of achieving better TWTT values, it comes at the cost of increased CPU times for the algorithm. This can be rationalized by considering that each execution of ALNS must handle a larger number of to-be-scheduled tasks, with lower values of $f$. Specifically, for the tight instances, increasing the value of $f$ from 0 to $\max (t^{\max},p^{\max})$ results in a $39.61\%$ reduction in solution time. When the value of $f$ is further increased to $t^{\max}+p^{\max}$ from $\max (t^{\max},p^{\max})$, the solution time decreases by $23.63\%$. This pattern remains consistent for the loose super-instances as well. In this case, increasing the value of $f$ from 0 to $\max (t^{\max},p^{\max})$ reduces the solution time by $16.63\%$, while increasing it to $t^{\max}+p^{\max}$ from $\max (t^{\max},p^{\max})$ results in a $1.27\%$ decrease in solution time.

\section{Conclusion}\label{section:conclusion}

In the realm of dynamic on-demand on-site service systems, optimizing operations in real-time is pivotal for businesses aiming to enhance their efficiency and productivity. In this article, we have embarked on a comprehensive exploration of various facets of a dynamic multi-skill workforce scheduling and routing problem with time windows and synchronization constraints, delving into reoptimization strategies and frozen period lengths. To address the outlined DWSRP-TW-SC, we developed an optimization framework that is triggered whenever a predetermined number of new tasks arrive. The initial step of this framework entails identifying the frozen tasks and establishing the earliest feasible time and location for the assigned team. Subsequently, the framework engages in
a re-optimization process for the ensuing Team Task Plan, aimed at minimizing the cumulative weighted completion time for all tasks. For the route redesign phase of this framework, we proposed two alternative methodologies: a mathematical model and a heuristic algorithm.

Our analysis began by assessing the effectiveness of the proposed ALNS compared to traditional solvers and a constructive heuristics. The ALNS algorithm demonstrated remarkable prowess in rapidly delivering high-quality solutions across a spectrum of DWSRP-TW-SC instances. It consistently outperformed the mathematical model and the constructive heuristic in terms of solution quality, often achieving optimal results within significantly shorter time frames. 

We further examined the impact of reoptimization strategies, discovering that a reduced frequency of reoptimization can lead to improved objective function values. While this improvement came at the cost of increased computation time, the trade-off was well-justified in scenarios where solution quality is of paramount importance.

Our investigation into the optimal frozen period length illuminated an intriguing balance between flexibility and stability. Extending the frozen period proved to be detrimental to solution quality, with shorter frozen periods resulting in significantly improved  the weighted sum of task throughput times. Decision-makers considering technology investments in reducing frozen periods should carefully assess the potential benefits, as even modest reductions can yield substantial enhancements.

In conclusion, our study provides insights for businesses dealing with dynamic workforce scheduling and routing problems. The proposed ALNS algorithm emerges as a powerful tool for real-time optimization, offering a compelling balance between solution quality and computation time. Furthermore, reoptimization strategies and frozen period lengths should all be tailored to suit specific business requirements, with careful consideration of the trade-offs involved. As industries continue to grapple with dynamic scheduling challenges, the principles and findings presented here can serve as a valuable compass guiding them toward more efficient and effective workforce management strategies. The interplay of dynamism, effective dynamism, and synchronization factors adds complexity to the dynamic on-demand on-site service systems, highlighting the need for holistic and adaptable solutions.

\section*{Acknowledgements}
This research was supported by TUBITAK [grant number 117M577].

During the preparation of this work the author(s) used ChatGBT in order to improve the language and readability. After using this tool/Service, the authors reveiewed and edited the content as needed and take full reponsibility for the content of the publication.
\newpage
\bibliography{ms}

\begin{thebibliography}{42}
\expandafter\ifx\csname natexlab\endcsname\relax\def\natexlab#1{#1}\fi
\providecommand{\url}[1]{\texttt{#1}}
\providecommand{\href}[2]{#2}
\providecommand{\path}[1]{#1}
\providecommand{\DOIprefix}{doi:}
\providecommand{\ArXivprefix}{arXiv:}
\providecommand{\URLprefix}{URL: }
\providecommand{\Pubmedprefix}{pmid:}
\providecommand{\doi}[1]{\href{http://dx.doi.org/#1}{\path{#1}}}
\providecommand{\Pubmed}[1]{\href{pmid:#1}{\path{#1}}}
\providecommand{\bibinfo}[2]{#2}
\ifx\xfnm\relax \def\xfnm[#1]{\unskip,\space#1}\fi
\bibitem[{Attanasio et~al.(2004)Attanasio, Cordeau, Ghiani and Laporte}]{Attanasio2004}
\bibinfo{author}{Attanasio, A.}, \bibinfo{author}{Cordeau, J.F.}, \bibinfo{author}{Ghiani, G.}, \bibinfo{author}{Laporte, G.}, \bibinfo{year}{2004}.
\newblock \bibinfo{title}{Parallel tabu search heuristics for the dynamic multi-vehicle dial-a-ride problem}.
\newblock \bibinfo{journal}{Parallel Computing} \bibinfo{volume}{30}, \bibinfo{pages}{377--387}.
\newblock \DOIprefix\doi{10.1016/j.parco.2003.12.001}.
\bibitem[{Barcelo et~al.(2007)Barcelo, Grzybowska and Pardo}]{Barcelo2007}
\bibinfo{author}{Barcelo, J.}, \bibinfo{author}{Grzybowska, H.}, \bibinfo{author}{Pardo, S.}, \bibinfo{year}{2007}.
\newblock \bibinfo{title}{Vehicle Routing And Scheduling Models, Simulation And City Logistics}. volume~\bibinfo{volume}{38}.
\newblock pp. \bibinfo{pages}{163--195}.
\newblock \DOIprefix\doi{10.1007/978-0-387-71722-7_8}.
\bibitem[{Bent and Van~Hentenryck(2004)}]{Bent2004}
\bibinfo{author}{Bent, R.}, \bibinfo{author}{Van~Hentenryck, P.}, \bibinfo{year}{2004}.
\newblock \bibinfo{title}{Scenario-based planning for partially dynamic vehicle routing with stochastic customers}.
\newblock \bibinfo{journal}{Operations Research} \bibinfo{volume}{52}, \bibinfo{pages}{977--987}.
\newblock \DOIprefix\doi{10.1287/opre.1040.0124}.
\bibitem[{{Benyahia} and {Potvin}(1998)}]{Benyahia1998}
\bibinfo{author}{{Benyahia}, I.}, \bibinfo{author}{{Potvin}, J..}, \bibinfo{year}{1998}.
\newblock \bibinfo{title}{Decision support for vehicle dispatching using genetic programming}.
\newblock \bibinfo{journal}{IEEE Transactions on Systems, Man, and Cybernetics - Part A: Systems and Humans} \bibinfo{volume}{28}, \bibinfo{pages}{306--314}.
\bibitem[{Berbeglia et~al.(2010)Berbeglia, Cordeau and Laporte}]{BERBEGLIA20108}
\bibinfo{author}{Berbeglia, G.}, \bibinfo{author}{Cordeau, J.F.}, \bibinfo{author}{Laporte, G.}, \bibinfo{year}{2010}.
\newblock \bibinfo{title}{Dynamic pickup and delivery problems}.
\newblock \bibinfo{journal}{European Journal of Operational Research} \bibinfo{volume}{202}, \bibinfo{pages}{8 -- 15}.
\newblock \URLprefix \url{http://www.sciencedirect.com/science/article/pii/S0377221709002999}, \DOIprefix\doi{https://doi.org/10.1016/j.ejor.2009.04.024}.
\bibitem[{Bieding et~al.(2009)Bieding, Görtz and Klose}]{Bieding2009}
\bibinfo{author}{Bieding, T.}, \bibinfo{author}{Görtz, S.}, \bibinfo{author}{Klose, A.}, \bibinfo{year}{2009}.
\newblock \bibinfo{title}{On line Routing per Mobile Phone A Case on Subsequent Deliveries of Newspapers}. volume \bibinfo{volume}{619}.
\newblock pp. \bibinfo{pages}{29--51}.
\newblock \DOIprefix\doi{10.1007/978-3-540-92944-4_3}.
\bibitem[{Borenstein et~al.(2008)Borenstein, Shah, Tsang, Dorne, Alsheddy and Voudouris}]{Borenstein2008}
\bibinfo{author}{Borenstein, Y.}, \bibinfo{author}{Shah, N.}, \bibinfo{author}{Tsang, E.}, \bibinfo{author}{Dorne, R.}, \bibinfo{author}{Alsheddy, A.}, \bibinfo{author}{Voudouris, C.}, \bibinfo{year}{2008}.
\newblock \bibinfo{title}{On the partitioning of dynamic scheduling problems -}, pp. \bibinfo{pages}{1691--1692}.
\newblock \DOIprefix\doi{10.1145/1389095.1389412}.
\bibitem[{Bostel et~al.(2008)Bostel, Dejax, Guez and Tricoire}]{bostel2008multiperiod}
\bibinfo{author}{Bostel, N.}, \bibinfo{author}{Dejax, P.}, \bibinfo{author}{Guez, P.}, \bibinfo{author}{Tricoire, F.}, \bibinfo{year}{2008}.
\newblock \bibinfo{title}{Multiperiod planning and routing on a rolling horizon for field force optimization logistics}, in: \bibinfo{booktitle}{The vehicle routing problem: latest advances and new challenges}. \bibinfo{publisher}{Springer}, pp. \bibinfo{pages}{503--525}.
\bibitem[{Campbell and Savelsbergh(2005)}]{Campbell2005}
\bibinfo{author}{Campbell, A.}, \bibinfo{author}{Savelsbergh, M.}, \bibinfo{year}{2005}.
\newblock \bibinfo{title}{Decision support for consumer direct grocery initiatives}.
\newblock \bibinfo{journal}{Transportation Science} \bibinfo{volume}{39}, \bibinfo{pages}{313--327}.
\newblock \DOIprefix\doi{10.1287/trsc.1040.0105}.
\bibitem[{Chang et~al.(2003)Chang, Hsueh and Student}]{Chang2003}
\bibinfo{author}{Chang, M.S.}, \bibinfo{author}{Hsueh, C.F.}, \bibinfo{author}{Student, P.}, \bibinfo{year}{2003}.
\newblock \bibinfo{title}{Real-time vehicle routing problem with time windows and simultaneous delivery/pickup demands}.
\newblock \bibinfo{journal}{Journal of the Eastern Asia Society for Transportation Studies} \bibinfo{volume}{5}.
\bibitem[{Cheung et~al.(2008)Cheung, Choy, Li, Shi and Tang}]{Cheung2008}
\bibinfo{author}{Cheung, B.}, \bibinfo{author}{Choy, K.}, \bibinfo{author}{Li, C.L.}, \bibinfo{author}{Shi, W.}, \bibinfo{author}{Tang, J.}, \bibinfo{year}{2008}.
\newblock \bibinfo{title}{Dynamic routing model and solution methods for fleet management with mobile technologies}.
\newblock \bibinfo{journal}{International Journal of Production Economics} \bibinfo{volume}{113}, \bibinfo{pages}{694--705}.
\newblock \DOIprefix\doi{10.1016/j.ijpe.2007.10.018}.
\bibitem[{Colby and Bell(2017)}]{colby_bell_2017}
\bibinfo{author}{Colby, C.}, \bibinfo{author}{Bell, K.}, \bibinfo{year}{2017}.
\newblock \bibinfo{title}{The on-demand economy is growing, and not just for the young and wealthy}.
\newblock \URLprefix \url{https://hbr.org/2016/04/the-on-demand-economy-is-growing-and-not-just-for-the-young-and-wealthy}.
\bibitem[{Cordeau et~al.(2007)Cordeau, Laporte, Potvin and Savelsbergh}]{CordeauLaporte2007}
\bibinfo{author}{Cordeau, J.}, \bibinfo{author}{Laporte, G.}, \bibinfo{author}{Potvin, J.}, \bibinfo{author}{Savelsbergh, M.}, \bibinfo{year}{2007}.
\newblock \bibinfo{title}{Transportation on Demand}. volume~\bibinfo{volume}{14}.
\newblock pp. \bibinfo{pages}{429 -- 466}.
\bibitem[{Cordeau and Laporte(2007)}]{CordeauLaporte2007_darp}
\bibinfo{author}{Cordeau, J.F.}, \bibinfo{author}{Laporte, G.}, \bibinfo{year}{2007}.
\newblock \bibinfo{title}{The dial-a-ride problem (darp): Models and algorithms}.
\newblock \bibinfo{journal}{Annals OR} \bibinfo{volume}{153}, \bibinfo{pages}{29--46}.
\newblock \DOIprefix\doi{10.1007/s10479-007-0170-8}.
\bibitem[{Coslovich et~al.(2006)Coslovich, Pesenti and Ukovich}]{COSLOVICH20061605}
\bibinfo{author}{Coslovich, L.}, \bibinfo{author}{Pesenti, R.}, \bibinfo{author}{Ukovich, W.}, \bibinfo{year}{2006}.
\newblock \bibinfo{title}{A two-phase insertion technique of unexpected customers for a dynamic dial-a-ride problem}.
\newblock \bibinfo{journal}{European Journal of Operational Research} \bibinfo{volume}{175}, \bibinfo{pages}{1605 -- 1615}.
\newblock \URLprefix \url{http://www.sciencedirect.com/science/article/pii/S0377221705002262}, \DOIprefix\doi{https://doi.org/10.1016/j.ejor.2005.02.038}.
\bibitem[{Demiray et~al.(2023)Demiray, Yucel and Tolga}]{mendeley_data}
\bibinfo{author}{Demiray, O.}, \bibinfo{author}{Yucel, E.}, \bibinfo{author}{Tolga, D.}, \bibinfo{year}{2023}.
\newblock \bibinfo{title}{Dwsrp data sets, mendeley data}.
\newblock \URLprefix \url{https://data.mendeley.com/datasets/rg9pshry4m/1}.
\bibitem[{Ferrucci and Bock(2014)}]{Ferrucci2014}
\bibinfo{author}{Ferrucci, F.}, \bibinfo{author}{Bock, S.}, \bibinfo{year}{2014}.
\newblock \bibinfo{title}{Real-time control of express pickup and delivery processes in a dynamic environment}.
\newblock \bibinfo{journal}{Transportation Research Part B: Methodological} \bibinfo{volume}{63}, \bibinfo{pages}{1–14}.
\newblock \DOIprefix\doi{10.1016/j.trb.2014.02.001}.
\bibitem[{Gendreau et~al.(2006)Gendreau, Guertin, Potvin and Séguin}]{Gendreau2006}
\bibinfo{author}{Gendreau, M.}, \bibinfo{author}{Guertin, F.}, \bibinfo{author}{Potvin, J.Y.}, \bibinfo{author}{Séguin, R.}, \bibinfo{year}{2006}.
\newblock \bibinfo{title}{Neighborhood search heuristics for a dynamic vehicle dispatching problem with pick-ups and deliveries}.
\newblock \bibinfo{journal}{Transportation Research Part C: Emerging Technologies} \bibinfo{volume}{14}, \bibinfo{pages}{157--174}.
\newblock \DOIprefix\doi{10.1016/j.trc.2006.03.002}.
\bibitem[{Gendreau et~al.(1999)Gendreau, Guertin, Potvin and Taillard}]{Gendreau1999}
\bibinfo{author}{Gendreau, M.}, \bibinfo{author}{Guertin, F.}, \bibinfo{author}{Potvin, J.Y.}, \bibinfo{author}{Taillard, E.}, \bibinfo{year}{1999}.
\newblock \bibinfo{title}{Parallel tabu search for real-time vehicle routing and dispatching}.
\newblock \bibinfo{journal}{Transportation Science} \bibinfo{volume}{33}, \bibinfo{pages}{381--390}.
\newblock \DOIprefix\doi{10.1287/trsc.33.4.381}.
\bibitem[{Haghani and Jung(2005)}]{HAGHANI2005}
\bibinfo{author}{Haghani, A.}, \bibinfo{author}{Jung, S.}, \bibinfo{year}{2005}.
\newblock \bibinfo{title}{A dynamic vehicle routing problem with time-dependent travel times}.
\newblock \bibinfo{journal}{Computers \& Operations Research} \bibinfo{volume}{32}, \bibinfo{pages}{2959 -- 2986}.
\newblock \URLprefix \url{http://www.sciencedirect.com/science/article/pii/S0305054804000887}, \DOIprefix\doi{https://doi.org/10.1016/j.cor.2004.04.013}.
\bibitem[{van Hemert and la(2004)}]{Hemert2004}
\bibinfo{author}{van Hemert, J.}, \bibinfo{author}{la, J.}, \bibinfo{year}{2004}.
\newblock \bibinfo{title}{Dynamic routing problems with fruitful regions: Models and evolutionary computation}.
\newblock \DOIprefix\doi{10.1007/978-3-540-30217-9_70}.
\bibitem[{Horn(2002)}]{Horn2002}
\bibinfo{author}{Horn, M.}, \bibinfo{year}{2002}.
\newblock \bibinfo{title}{Multi-modal and demand-responsive passenger transport systems: A modelling framework with embedded control systems}.
\newblock \bibinfo{journal}{Transportation Research Part A: Policy and Practice} \bibinfo{volume}{36}, \bibinfo{pages}{167--188}.
\newblock \DOIprefix\doi{10.1016/S0965-8564(00)00043-4}.
\bibitem[{Ichoua et~al.(2000)Ichoua, Gendreau and Potvin}]{Ichoua2000}
\bibinfo{author}{Ichoua, S.}, \bibinfo{author}{Gendreau, M.}, \bibinfo{author}{Potvin, J.Y.}, \bibinfo{year}{2000}.
\newblock \bibinfo{title}{Diversion issues in real-time vehicle dispatching}.
\newblock \bibinfo{journal}{Transportation Science} \bibinfo{volume}{34}, \bibinfo{pages}{426--438}.
\newblock \DOIprefix\doi{10.1287/trsc.34.4.426.12325}.
\bibitem[{Ichoua et~al.(2003)Ichoua, Gendreau and Potvin}]{Ichoua2003}
\bibinfo{author}{Ichoua, S.}, \bibinfo{author}{Gendreau, M.}, \bibinfo{author}{Potvin, J.Y.}, \bibinfo{year}{2003}.
\newblock \bibinfo{title}{Vehicle dispatching with time - dependent travel times}.
\newblock \bibinfo{journal}{European Journal of Operational Research} \bibinfo{volume}{144}, \bibinfo{pages}{379--396}.
\newblock \DOIprefix\doi{10.1016/S0377-2217(02)00147-9}.
\bibitem[{Ichoua et~al.(2007)Ichoua, Gendreau and Potvin}]{Ichoua2007}
\bibinfo{author}{Ichoua, S.}, \bibinfo{author}{Gendreau, M.}, \bibinfo{author}{Potvin, J.Y.}, \bibinfo{year}{2007}.
\newblock \bibinfo{title}{Planned Route Optimization For Real-Time Vehicle Routing}. volume~\bibinfo{volume}{38}.
\newblock pp. \bibinfo{pages}{1--18}.
\newblock \DOIprefix\doi{10.1007/978-0-387-71722-7_1}.
\bibitem[{Larsen(2001)}]{Larsen2001}
\bibinfo{author}{Larsen, A.}, \bibinfo{year}{2001}.
\newblock \bibinfo{title}{The Dynamic Vehicle Routing Problem}.
\newblock Ph.D. thesis. Department of Mathematical Modelling (IMM), TechnicalUniversity of Denmark (DTU).
\bibitem[{Lois and Ziliaskopoulos(2017)}]{LOIS2017377}
\bibinfo{author}{Lois, A.}, \bibinfo{author}{Ziliaskopoulos, A.}, \bibinfo{year}{2017}.
\newblock \bibinfo{title}{Online algorithm for dynamic dial a ride problem and its metrics}.
\newblock \bibinfo{journal}{Transportation Research Procedia} \bibinfo{volume}{24}, \bibinfo{pages}{377 -- 384}.
\newblock \URLprefix \url{http://www.sciencedirect.com/science/article/pii/S2352146517303782}, \DOIprefix\doi{https://doi.org/10.1016/j.trpro.2017.05.097}. \bibinfo{note}{3rd Conference on Sustainable Urban Mobility, 3rd CSUM 2016, 26 – 27 May 2016, Volos, Greece}.
\bibitem[{Lund et~al.(1996)Lund, Madsen and J.M.}]{Lund1996}
\bibinfo{author}{Lund, K.}, \bibinfo{author}{Madsen, O.}, \bibinfo{author}{J.M., R.}, \bibinfo{year}{1996}.
\newblock \bibinfo{title}{Vehicle routing problems with varying degrees of dynamism} .
\bibitem[{Montemanni et~al.(2005)Montemanni, Gambardella, Rizzoli and Donati}]{Montemanni2005}
\bibinfo{author}{Montemanni, R.}, \bibinfo{author}{Gambardella, L.M.}, \bibinfo{author}{Rizzoli, A.E.}, \bibinfo{author}{Donati, A.}, \bibinfo{year}{2005}.
\newblock \bibinfo{title}{Ant colony system for a dynamic vehicle routing problem}.
\newblock \bibinfo{journal}{J. Comb. Optim.} \bibinfo{volume}{10}, \bibinfo{pages}{327--343}.
\newblock \DOIprefix\doi{10.1007/s10878-005-4922-6}.
\bibitem[{Pillac et~al.(2013)Pillac, Gendreau, Guéret and Medaglia}]{Pillac2013}
\bibinfo{author}{Pillac, V.}, \bibinfo{author}{Gendreau, M.}, \bibinfo{author}{Guéret, C.}, \bibinfo{author}{Medaglia, A.}, \bibinfo{year}{2013}.
\newblock \bibinfo{title}{A review of dynamic vehicle routing problems}.
\newblock \bibinfo{journal}{European Journal of Operational Research} \bibinfo{volume}{225}, \bibinfo{pages}{1–11}.
\newblock \DOIprefix\doi{10.1016/j.ejor.2012.08.015}.
\bibitem[{Pillac et~al.(2018)Pillac, Guéret and Medaglia}]{Pillac2018}
\bibinfo{author}{Pillac, V.}, \bibinfo{author}{Guéret, C.}, \bibinfo{author}{Medaglia, A.}, \bibinfo{year}{2018}.
\newblock \bibinfo{title}{A Fast Reoptimization Approach for the Dynamic Technician Routing and Scheduling Problem}.
\newblock pp. \bibinfo{pages}{347--367}.
\newblock \DOIprefix\doi{10.1007/978-3-319-58253-5_20}.
\bibitem[{Pisinger and Ropke(2007)}]{PisingerRopke2007}
\bibinfo{author}{Pisinger, D.}, \bibinfo{author}{Ropke, S.}, \bibinfo{year}{2007}.
\newblock \bibinfo{title}{A general heuristic for vehicle routing problems}.
\newblock \bibinfo{journal}{Computers \& Operations research} \bibinfo{volume}{34}, \bibinfo{pages}{2403--2435}.
\bibitem[{Potvin and Rousseau(1993)}]{PotvinRousseau1993}
\bibinfo{author}{Potvin, J.Y.}, \bibinfo{author}{Rousseau, J.M.}, \bibinfo{year}{1993}.
\newblock \bibinfo{title}{A parallel route building algorithm for the vehicle routing and scheduling problem with time windows}.
\newblock \bibinfo{journal}{European Journal of Operational Research} \bibinfo{volume}{66}, \bibinfo{pages}{331--340}.
\bibitem[{Regnier-Coudert et~al.(2016)Regnier-Coudert, Mccall, Ayodele and Anderson}]{Regnier-Coudert2016}
\bibinfo{author}{Regnier-Coudert, O.}, \bibinfo{author}{Mccall, J.}, \bibinfo{author}{Ayodele, M.}, \bibinfo{author}{Anderson, S.}, \bibinfo{year}{2016}.
\newblock \bibinfo{title}{Truck and trailer scheduling in a real world, dynamic and heterogeneous context}.
\newblock \bibinfo{journal}{Transportation Research Part E: Logistics and Transportation Review} \bibinfo{volume}{93}, \bibinfo{pages}{389--408}.
\newblock \DOIprefix\doi{10.1016/j.tre.2016.06.010}.
\bibitem[{Ropke and Pisinger(2006a)}]{RopkePisinger2006a}
\bibinfo{author}{Ropke, S.}, \bibinfo{author}{Pisinger, D.}, \bibinfo{year}{2006}a.
\newblock \bibinfo{title}{An adaptive large neighborhood search heuristic for the pickup and delivery problem with time windows}.
\newblock \bibinfo{journal}{Transportation Science} \bibinfo{volume}{40}, \bibinfo{pages}{455--472}.
\bibitem[{Ropke and Pisinger(2006b)}]{RopkePisinger2006b}
\bibinfo{author}{Ropke, S.}, \bibinfo{author}{Pisinger, D.}, \bibinfo{year}{2006}b.
\newblock \bibinfo{title}{A unified heuristic for a large class of vehicle routing problems with backhauls.}
\newblock \bibinfo{journal}{European Journal of Operational Research} \bibinfo{volume}{171}, \bibinfo{pages}{750--775}.
\bibitem[{Savelsbergh and Sol(1997)}]{Savelsbergh97}
\bibinfo{author}{Savelsbergh, M.}, \bibinfo{author}{Sol, M.}, \bibinfo{year}{1997}.
\newblock \bibinfo{title}{Drive: Dynamic routing of independent vehicles}.
\newblock \bibinfo{journal}{Operations Research} \bibinfo{volume}{46}.
\newblock \DOIprefix\doi{10.1287/opre.46.4.474}.
\bibitem[{Shaw(1997)}]{shaw1997}
\bibinfo{author}{Shaw, P.}, \bibinfo{year}{1997}.
\newblock \bibinfo{title}{A new local search algorithm providing high quality solutions to vehicle routing problems}.
\newblock \bibinfo{journal}{APES Group, Dept of Computer Science, University of Strathclyde, Glasgow, Scotland, UK} .
\bibitem[{Shaw(1998)}]{shaw1998}
\bibinfo{author}{Shaw, P.}, \bibinfo{year}{1998}.
\newblock \bibinfo{title}{Using constraint programming and local search methods to solve vehicle routing problems}, in: \bibinfo{booktitle}{International conference on principles and practice of constraint programming}, \bibinfo{organization}{Springer}. pp. \bibinfo{pages}{417--431}.
\bibitem[{Taillard et~al.(2001)Taillard, Gambardella, Gendreau and Potvin}]{Taillard2001}
\bibinfo{author}{Taillard, E.}, \bibinfo{author}{Gambardella, L.M.}, \bibinfo{author}{Gendreau, M.}, \bibinfo{author}{Potvin, J.Y.}, \bibinfo{year}{2001}.
\newblock \bibinfo{title}{Adaptive memory programming: A unified view of metaheuristics}.
\newblock \bibinfo{journal}{European Journal of Operational Research} \bibinfo{volume}{135}, \bibinfo{pages}{1--16}.
\newblock \DOIprefix\doi{10.1016/S0377-2217(00)00268-X}.
\bibitem[{Wilson and Colvin(1977)}]{WilsonColvin1977}
\bibinfo{author}{Wilson, N.H.M.}, \bibinfo{author}{Colvin, N.J.}, \bibinfo{year}{1977}.
\newblock \bibinfo{title}{Computer control of the rochester dial-a-ride system}.
\newblock \bibinfo{journal}{Cambridge: Massachusetts Institute of Technology, Center for Transportation Studies} .
\bibitem[{Çakırgil et~al.(2020)Çakırgil, Yücel and Kuyzu}]{Cakirgil2020}
\bibinfo{author}{Çakırgil, S.}, \bibinfo{author}{Yücel, E.}, \bibinfo{author}{Kuyzu, G.}, \bibinfo{year}{2020}.
\newblock \bibinfo{title}{An integrated solution approach for multi-objective, multi-skill workforce scheduling and routing problems}.
\newblock \bibinfo{journal}{Computers \& Operations Research} , \bibinfo{pages}{104908}\DOIprefix\doi{10.1016/j.cor.2020.104908}.

\end{thebibliography}

\end{document}